\documentclass[]{interact}

\usepackage[numbers,sort&compress]{natbib}
\bibpunct[, ]{[}{]}{,}{n}{,}{,}
\makeatletter
\def\NAT@def@citea{\def\@citea{\NAT@separator}}
\makeatother

\usepackage{mathtools}
\usepackage[bookmarks]{hyperref}
\usepackage{enumerate,longtable}
\usepackage{color,colortbl}
\usepackage{mathrsfs}
\usepackage{multirow}

\numberwithin{equation}{section}

\theoremstyle{plain}
\newtheorem{theorem}{Theorem}[section]

\theoremstyle{definition}

\newtheorem{example}[theorem]{Example}

\theoremstyle{remark}
\newtheorem{remark}[theorem]{Remark}

\newcommand{\A}{\mathcal{A}}
\newcommand{\Au}{\mathrm{Aut}}
\newcommand{\ad}{\mathrm{ad}}
\newcommand{\bal}{\begin{array}{l}}
\newcommand{\bc}{\begin{cases}}

\newcommand{\C}{\mathbb{C}}
\newcommand{\co}{\coloneqq}
\newcommand{\di}{\mathrm{diag}}
\newcommand{\ea}{\end{array}}
\newcommand{\ec}{\end{cases}}
\newcommand{\F}{\mathbb{F}}
\newcommand{\g}{\mathfrak{g}}

\newcommand{\h}{\mathfrak{h}}
\newcommand{\Ls}{\mathscr{L}}
\newcommand{\M}{\mathrm{Mat}}
\newcommand{\n}{\mathfrak{n}}

\newcommand{\R}{\mathbb{R}}
\newcommand{\ra}{\mathrm{rank}\,}
\newcommand{\Sp}{\mathrm{span}}

\newcommand{\xy}{[X, Y]}
\newcommand{\Z}{\mathcal{Z}}

\begin{document}


\title{Classification of 7-dimensional solvable Lie algebras having 5-dimensional nilradicals}

\author{
	\name{
		Vu A. Le\textsuperscript{a}, Tuan A. Nguyen\textsuperscript{b}, Tu T. C. Nguyen\textsuperscript{c,d}, 
		Tuyen T. M. Nguyen\textsuperscript{c,e}  and
		Thieu N. Vo\textsuperscript{f}\thanks{CONTACT Thieu N. Vo. Email: vongocthieu@tdtu.edu.vn}
	}
	\affil{
		\textsuperscript{a}Department of Economic Mathematics, University of Economics and Law, 
			Vietnam National University -- Ho Chi Minh City, Vietnam;
		\textsuperscript{b}Faculty of Political Science and Pedagogy, 
			Ho Chi Minh City University of Physical Education and Sport, Vietnam; 
		\textsuperscript{c}Faculty of Mathematics and Computer Science, 
			University of Science -- Vietnam National University, Ho Chi Minh City, Vietnam;
		\textsuperscript{d}College of Natural Sciences, Can Tho University, Can Tho City, Vietnam;
		\textsuperscript{e}Faculty of Mathematics and Computer Science Teacher Education, 
			Dong Thap University, Dong Thap Province, Vietnam;
		\textsuperscript{f}Fractional Calculus, Optimization and Algebra Research Group, 
			Faculty of Mathematics and Statistics, Ton Duc Thang University, Ho Chi Minh City, Vietnam
	}
}

\maketitle

\begin{abstract}
This paper presents a classification of 7-dimensional real and complex indecomposable 
solvable Lie algebras having some 5-dimensional nilradicals. Afterwards, we combine our 
results with those of Rubin and Winternitz (1993), Ndogmo and Winternitz (1994), 
\v Snobl and Winternitz (2005, 2009), \v Snobl and Kar\'asek (2010) to obtain a complete classification of 
7-dimensional real and complex indecomposable solvable Lie algebras with 5-dimensional nilradicals.
In association with Gong (1998), Parry (2007), Hindeleh and Thompson (2008), 
we achieve a classification of 7-dimensional real and complex indecomposable solvable Lie algebras. 
\end{abstract}


\begin{amscode}
	15A21; 16G60; 17B30; 20G05
\end{amscode}

\section{Introduction}\label{sec1}

Classifying Lie algebras is a central problem in Lie theory. By the well-known theorems of
Levi \cite{Lev05} and Malcev \cite{Mal45}, the problem of classifying Lie algebras over a field of characteristic zero 
is reduced to the problem of classifying semi-simple Lie algebras and solvable ones. 
Semi-simple Lie algebras were fully classified by Cartan~\cite{Car94} (over the complex field $\C$) 
and Ganmatcher \cite{Gan39} (over the real field $\R$). However, classifying solvable 
Lie algebras is much harder, and in general, it still remains open.

In history, the classification of complex and real solvable Lie algebras was achieved up to dimension 6 by Lie \cite{LE93}, 
Bianchi \cite{Bia03}, Dixmier \cite{Dix58}, Morozov \cite{Mor58}, Mubarakzyanov \cite{Mub63a,Mub63b,Mub63c} 
and Turkowski \cite{Tur90}. All of these results were summarized by \v Snobl and Winternitz \cite{SW14}. 
For complex and real solvable Lie algebras of dimensions higher than 6, several partial results are presented in
\cite{BEN94,BFN03,Che12,Gon98,KT90,Par07,HT08,See94,TK91,Tsa99,TKK00} and references therein. 
Most of these results concern with nilpotent Lie algebras. 
So far, a full classification of solvable Lie algebras of dimension 7 has not been completed yet.

In this paper, we present a full classification of 7-dimensional real and complex indecomposable solvable Lie algebras.
Our method is based on the fact that a given solvable Lie algebra $L$ can be considered as an extension 
of its nilradical $N(L)$, that is, the maximal nilpotent ideal of $L$. Therefore, we start from a nilpotent Lie algebra 
and classify all 7-dimensional solvable Lie algebras which admit it as their nilradical. 
This method, perhaps, was initialized in 1963 by a series of articles of Mubarakzyanov \cite{Mub63a,Mub63b,Mub63c} 
when he classified solvable Lie algebras of dimensions 4 and 5 over a field of characteristic zero. 
By using the same method, the results for the case of dimension 6 were also achieved by 
Mubarakzyanov \cite{Mub63c} and Turkowski \cite{Tur90}. Furthermore, results in \cite{NW94,RW93,ST13,Sha16,Sno11,SK10,SW05,SW09,SW12,WLD08} show that this method seems to be very effective.

In case $\dim L=7$, it follows from Mubarakzyanov \cite[Theorem 5]{Mub63a} that $\dim N(L) \in \{4, 5, 6, 7\}$. 
Therefore, the problem of classifying 7-dimensional solvable Lie algebras consists of four cases 
according to the possible values of $\dim N(L)$. Moreover, three cases in which $\dim N(L) \in \{4, 6, 7\}$ 
were considered by Hindeleh and Thompson \cite{HT08}, Parry \cite{Par07} and Gong \cite{Gon98}, respectively. 
To our knowledge, the remaining case when $\dim N(L)=5$ has just been solved partially, up to date.

According to Dixmier \cite[Proposition 1]{Dix58}, the class of 5-dimensional real and complex nilpotent 
Lie algebras consists of nine Lie algebras as follows:
\begin{itemize}
	\item $(\g_1)^5$: the 5-dimensional abelian Lie algebra;
	\item $(\g_1)^2 \oplus \g_3$: $[x_1, x_2] = x_3$;
	\item $\g_1 \oplus \g_4$: $[x_1, x_2] = x_3$, $[x_1, x_3] = x_4$;
	\item $\g_{5,1}$: $[x_1, x_2] = x_5$, $[x_3, x_4] = x_5$;
	\item $\g_{5,2}$: $[x_1, x_2] = x_4$, $[x_1, x_3] = x_5$;
	\item $\g_{5,3}$: $[x_1, x_2] = x_4$, $[x_1, x_4] = x_5$, $[x_2, x_3] = x_5$;
	\item $\g_{5,4}$: $[x_1, x_2] = x_3$, $[x_1, x_3] = x_4$, $[x_2, x_3] = x_5$;
	\item $\g_{5,5}$: $[x_1, x_2] = x_3$, $[x_1, x_3] = x_4$, $[x_1, x_4] = x_5$;
	\item $\g_{5,6}$: $[x_1, x_2] = x_3$, $[x_1, x_3] = x_4$, $[x_1, x_4] = x_5$, $[x_2, x_3] = x_5$.
\end{itemize}
Finite-dimensional real and complex indecomposable solvable extensions of $(\g_1)^5$, $\g_{5,1}$, $\g_{5,3}$, $\g_{5,5}$ 
and $\g_{5,6}$ were considered by Ndogmo and Winternitz \cite{NW94}, Rubin and Winternitz \cite{RW93}, 
\v Snobl and Kar\'asek \cite{SK10}, \v Snobl and Winternitz \cite{SW05,SW09}, respectively. Finite-dimensional 
complex indecomposable solvable extensions of $\g_1 \oplus \g_4$ were considered by Wang et al. \cite{WLD08}.

The main goal of this paper is to classify all 7-dimensional indecomposable solvable extensions 
of $(\g_1)^2 \oplus \g_3$, $\g_1 \oplus \g_4$, $\g_{5,2}$ and $\g_{5,4}$ (see Theorem \ref{thm4.1}). 
Then, we combine our results with those of \cite{NW94,RW93,SK10,SW05,SW09} to obtain a classification 
of 7-dimensional indecomposable Lie algebras having 5-dimensional nilradicals (see Theorem \ref{thm4.2}). 
In association with Gong \cite{Gon98}, Parry \cite{Par07}, Hindeleh and Thompson \cite{HT08}, 
we achieve a full classification of 7-dimensional indecomposable solvable Lie algebras.

We organize the paper into five sections and one appendix. Section \ref{sec2} describes the classification procedure.
In Section \ref{sec3}, we give explicit computations for a sample case of nilradical. 
Afterwards, we formulate two main theorems of the paper in Section \ref{sec4}. 
Section \ref{sec5} contains a full summary for the classification of 7-dimensional solvable Lie algebras. 
Finally, the appendix presents the full lists of Lie algebras achieved in our classification with precisely isomorphic conditions.

\section{The procedure of classification}\label{sec2}

The classification in this paper proceeds in two stages which will be described immediately as follows.
From now on, $\F$ will be $\R$ or $\C$.

\subsection{Construction of Lie algebras}\label{ssec2.1}

In the first stage, we construct four lists $\Ls_1$, $\Ls_2$, $\Ls_3$ and $\Ls_4$ which consist 
of 7-dimensional indecomposable solvable $\F$-Lie algebras $L$ having nilradicals 
$(\g_1)^2 \oplus \g_3$, $\g_1 \oplus \g_4$, $\g_{5,2}$ and $\g_{5,4}$, respectively. 
Note that the solvable extension of a given nilpotent Lie algebra is standard 
and can be found in many textbooks (see, e.g., \cite{SW14}).

First of all, we fix a Lie algebra in $\left\lbrace (\g_1)^2 \oplus \g_3, \g_1 \oplus \g_4, \g_{5,2}, \g_{5,4} \right\rbrace$ 
which plays a role as the input nilradical $N(L)$. The basis of $N(L)$ is always assume to be $\{X_1, \ldots, X_5\}$. 
By adding to the basis $\{X_1, \ldots, X_5\}$ two elements, say $X$ and $Y$, we obtain a basis $\{X_1, \ldots, X_5, X, Y\}$ 
of $L$. Then, Lie brackets of $L$ are absolutely determined by $\xy$, $[X, X_i]$ and $[Y, X_i]$ for $i = 1, \ldots, 5$.
Since the derived algebra of a solvable Lie algebra is contained in its nilradical (see \cite[Chapter II, Section 7, Corollary 1]{Jac62}), 
these Lie brackets can be represented as follows:
\[
	\begin{array}{l l l l}
		\xy = \sum \limits_{j=1}^5 \sigma_jX_j, &
		\left[X, X_i\right] = \sum \limits_{j=1}^5 a_{ij}X_j, &
		\left[Y, X_i\right] = \sum \limits_{j=1}^5 b_{ij}X_j; & 1 \leq i \leq 5.
	\end{array}
\]
Set $A \co \left(a_{ij}\right)$ and $B \co \left(b_{ij}\right)$. We call $A, B \in \M_5(\F)$
the \emph{structure matrices} of $L$. Then, all we have to do is to determine all possibilities 
of structure constants $\sigma_i \in \F$ and the structure matrices $A, B$. To this end, the following techniques will be used.
\begin{enumerate}
	\item First of all, twenty Jacobi identities involving $(X, X_i, X_j)$ and $(Y, X_i, X_j)$ initialize the original forms 
	of $A$ and $B$, respectively.
	
	\item Next, five Jacobi identities involving $(X, Y, X_i)$ construct a relation between $A$ and $B$. 
	Moreover, $[A, B]$ is an inner derivation of $N(L)$, i.e., we have
	\[
		\begin{array}{l l}
			[A, B] = \sum \limits_{i=1}^5 \sigma_i a_{X_i}, & a_{X_i} \co \left(\ad_{X_i}\right)^T\vert_{N(L)},
		\end{array}
	\]
	where $\ad$ is the adjoint operator and the superscript $T$ indicates the transpose of a matrix.
		
	\item Two supplemented elements $X$ and $Y$ must be \emph{linearly nil-independent} to ensure 
	that the dimension of $N(L)$ is not larger than 5. This is equivalent to the fact that $A$ and $B$ are 
	linearly nil-independent, i.e., if $\alpha A + \beta B$ is a nilpotent matrix then $\alpha = \beta = 0$.
	
	\item We use the three following types of transformations alternatively not only to eliminate $\sigma_i$ or normalize 
	$\sigma_i \neq 0$ but also to simplify $A$ and $B$:
	\begin{itemize}
		\item Modifying $A$ and $B$ by inner derivations of $N(L)$, i.e., $A' = A + \sum \limits_{i=1}^5 \alpha_ia_{X_i}$, 
		$B' = B + \sum \limits_{i=1}^5 \beta_ia_{X_i}$ with $\alpha_i, \beta_i \in \F$.
		
		\item The automorphisms of $N(L)$: if $G \in \Au(N(L))$ then it will transforms the structure matrices into 
		$A' = GAG^{-1}$ and $B' = GBG^{-1}$.
		
		\item The last one is the transformation concerning with $X$ and $Y$.
	\end{itemize}
	
	\item Besides, due to Mubarakzyanov \cite[Corollary 2]{Mub66}, the following relation holds for finite-dimensional 
	solvable Lie algebras $L$ over a field of characteristic zero
	\begin{equation}\label{Mub-ineq2}
		2\dim N(L) \geq \dim L + \dim \Z,
	\end{equation}
	where $\Z$ is the center of $L$. Inequation \eqref{Mub-ineq2} gives an upper bound of $\dim \Z$ 
	which is very useful to decide which cases can happen.
\end{enumerate}

Afterwards, we repeat all techniques above, case by case, for other nilradicals in 
$\left\lbrace (\g_1)^2 \oplus \g_3, \g_1 \oplus \g_4, \g_{5,2}, \g_{5,4}\right\rbrace$. 
By this way, we can obtain four lists $\Ls_i$ as desired.

\subsection{Testing isomorphism for the obtained Lie algebras}\label{ssec2.2}

To optimize four lists $\Ls_i$, we need to test isomorphism between Lie algebras in the lists 
as well as refine parameters, if any. This stage is necessary since it makes $\Ls_i$ more compact. 
Moreover, we can avoid redundancy, i.e., different Lie algebras in $\Ls_i$ are non-isomorphic. 

First of all, Lie algebras with different invariants (such as the dimensions of ideals in characteristic series 
or the dimensions of centers) are non-isomorphic. For Lie algebras sharing same invariants, 
we use computer algebra tools to verify their isomorphism. Gerdt and Lassner \cite{GL93}, perhaps, 
were the first authors considering the problem of testing isomorphism of real and complex Lie algebras 
from a computer algebra point of view. They reduce the problem of testing isomorphism of Lie algebras 
to the problem of testing the existence of solutions of a polynomial system. Gr\"obner basis technique is 
then used to solve the latter problem. However, since the complexity of computing Gr\"obner bases is very costly, 
this method is impractical when the dimension pass 6, especially, in case of parametric Lie algebras.

In this paper, we use another computer algebra tool which is the so-called \emph{triangular decomposition} 
instead of Gr\"obner bases. Following the idea of Gerdt and Lassner \cite{GL93}, we also re-write isomorphic 
conditions between Lie algebras, even if parametric Lie algebras, in terms of polynomial systems and 
semi-algebraic systems. Afterwards, we use triangular decomposition to decide whether these systems admit roots or not; 
and if they do, we can find explicit roots to construct isomorphisms. This testing isomorphism procedure is specified by 
algorithms that are run by Maple software with supports of a hyper-computer. Details of these algorithms can be found in \cite{NLV21}.
We also note that these algorithms in \cite{NLV21} are valid over fields of characteristic not only 0 but also prime.
In this paper, we simply use them over $\F$.

\begin{remark}\label{rem2.1}
	This stage has a further advantage as follows. Assume that there is a list $\Ls$ which consists of 
	\emph{real} Lie algebras	satisfying certain properties. Since each algebra in $\Ls$ can be seen as a complex one, 
	we first sweep out algebras that do not exist over $\C$ by considering structure matrices. 
	Afterwards, we check indecomposability since an indecomposable algebra over $\R$ may be decomposable over $\C$. 
	Finally, we test their isomorphism over $\C$. Consequently, a similar list over $\C$ can be derived from $\Ls$.
	Below is a demonstrative example.
\end{remark}

\begin{example}\label{ex2.2}
	Parry \cite{Par07} in 2007 classified 7-dimensional real indecomposable solvable
	Lie algebras with codimension one nilradicals. In case of nilradical $\R \oplus \g_{5,1}$, there is a
	family, namely $[7, [6, 5], 1, 3]$, which is as follows (see \cite[Appendix B.6]{Par07}):
	$[e_3, e_5] = e_2$, $[e_4, e_6] = e_2$, $[e_1, e_7] = ae_1$, $[e_3, e_7] = -e_3$,
	$[e_4, e_7] = -ae_4$, $[e_5, e_7] = e_5$, $[e_6, e_7] = -e_1 + ae_6$ with $a \in \R$.
	We claim that the real parameter $a$ can be reduced to $a \geq 0$ instead of $a \in \R$. 
	In fact, the testing isomorphism procedure above shows that all algebras depending on $a$ 
	are isomorphic to those on $-a$ by the following isomorphism (we omit off-diagonal zeros):
	\[
	 	\begin{bmatrix}
	 		-1 \\
	 		& 1 \\
	 		&& 0 && 1 \\ 
	 		&&& 1 \\
	 		&& -1 && 0 \\
	 		&&&&& 1 \\
	 		&&&&&& -1
	 	\end{bmatrix}.
	\]
	Therefore, we can remove $[7, [6, 5], 1, 3]_{a < 0}$ to avoid redundancy. We have not checked all other algebras 
	in Parry's list, however, similar situations may occur, especially, for algebras containing more than one parameter. 
	Furthermore, if we consider $[7, [6, 5], 1, 3]$ over $\C$, the similar thing also happens: $[7, [6, 5], 1, 3]$ is indecomposable, 
	$a$ and $-a$ also determine isomorphic complex algebras, and we thus can reduce the complex parameter 
	$a$ to $a =0$ or $a \neq 0$ with $0 < \mathrm{arg}(a) \leq \pi$.
\end{example}

\section{A sample case: solvable extension of $\g_{5,2}$}\label{sec3}

The goal of this section is to classify all 7-dimensional indecomposable solvable $\F$-Lie algebras 
having niradical $\g_{5,2}$ with detailed computations by the procedure pointed out in Section \ref{sec2}. 
Recall that $\g_{5,2}$ is in Section \ref{sec1}.

First, twenty Jacobi identities involving $(X,X_i,X_j)$, $(Y,X_i,X_j)$ and transformations 
$X' \co X-a_{24}X_1+a_{14}X_2+a_{15}X_3$, $Y' \co Y-b_{24}X_1+b_{14}X_2+b_{15}X_3$ give:
\[
	\begin{array}{l l l}
		A=\begin{bmatrix} a & d & e \\ & b & f && g \\ & h & c & k & l \\ &&& a+b & f \\ &&& h & a+c \end{bmatrix}, &
		B=\begin{bmatrix} u & p & q \\ & v & r && t \\ & x & w & y & z \\ &&& u+v & r \\ &&& x & u+w \end{bmatrix}.
	\end{array}
\]

Next, five Jacobi identities involving $(X,Y,X_i)$ give $\sigma_1 = 0$ and
\begin{equation}\label{3.1}
	\begin{cases}
		(a-b)p+ex=(u-v)d+hq\\
		(a-c)q+dr=(u-w)e+pf \\
		fx=hr \\ az=ul \\ (b-c)r=(v-w)f\\
		fy+gx=ht+kr\\
		(a-b+c)t+lr=(u-v+w)g+fz \\
		(b-c)x=(v-w)h\\
		(a+b-c)y+zh=(u+v-w)k+lx
	\end{cases}		
\end{equation}
Moreover, $[A, B] = -\sigma_2a_{X_2} - \sigma_3a_{X_3}$, where $a_{X_i} \co \left(\ad_{X_i}\right)^T\vert_{\g_{5,2}}$. 
Put $C = \begin{bmatrix} b & f \\ h & c \end{bmatrix}$ and $D=\begin{bmatrix} v & r \\ x & w \end{bmatrix}$.
Then the third, fifth and eighth equations of \eqref{3.1} imply that $C$ and $D$ commute. 
Therefore, we can choose $\alpha, \beta, \gamma, \delta$ in the following automorphism of $\g_{5,2}$
\[
	\begin{array}{l l l}
		\begin{bmatrix} X'_1 \\ X'_2 \\ X'_3 \\ X'_4 \\ X'_5 \end{bmatrix} = G\begin{bmatrix} X_1 \\ X_2 \\ X_3 \\ X_4 \\ X_5 \end{bmatrix},
		&& G = \begin{bmatrix} 1 \\ & \alpha & \beta \\ &\gamma & \delta \\ &&& \alpha & \beta \\ &&& \gamma & \delta \end{bmatrix}
			\in \Au(\g_{5,2}),
	\end{array}
\]
such that the pair $(C, D)$ can obtain the three following types:
\[
	\begin{array}{l l l}
		\left(\begin{bmatrix} \lambda_1 & 0 \\ 0 & \lambda_2 \end{bmatrix}, \begin{bmatrix} \mu_1 & 0 \\ 0 & \mu_2 \end{bmatrix}\right); &
		\left(\begin{bmatrix} \lambda & 1 \\ 0 & \lambda \end{bmatrix}, \begin{bmatrix} \mu & \mu_1 \\ 0 & \mu \end{bmatrix} \right); &
		\left(\begin{bmatrix} \lambda_1 & \lambda_2 \\ -\lambda_2 & \lambda_1 \end{bmatrix}, \begin{bmatrix} \mu_1 & \mu_2 \\ -\mu_2 & \mu_1\end{bmatrix} \right), \; \lambda_2 \neq  0.
	\end{array}
\]
We do not need to consider permutations of three types above since they will return to three original ones 
if we interchange $X \leftrightarrow Y$. Now, these three types lead to three forms of structure matrices as follows:
\begin{enumerate}[\bf (1)]
	\item\label{form1} $A = \begin{bmatrix} a & d & e \\ & \lambda_1 &&& g \\ && \lambda_2 & k & l \\ &&& a+\lambda_1 \\ &&&& a+\lambda_2 \end{bmatrix}$ and $B=\begin{bmatrix} u & p & q \\ & \mu_1 &&& t \\ && \mu_2 & y & z \\ &&& u+\mu_1 \\ &&&& u+\mu_2 \end{bmatrix}$;
	\item \label{form2} $A = \begin{bmatrix} a & d & e \\ & \lambda & 1 && g \\ && \lambda & k & l \\ &&& a+\lambda & 1 \\ &&&& a+\lambda \end{bmatrix}$ and $B = \begin{bmatrix} u & p & q \\ & \mu & \mu_1 && t \\ 	&& \mu & y & z \\ &&& u+\mu & \mu_1 \\ &&&& u+\mu \end{bmatrix}$;
	\item\label{form3} $A = \begin{bmatrix} a & d & e \\ & \lambda_1 & \lambda_2 && g \\ & -\lambda_2 & \lambda_1 & k & l \\ &&& a+\lambda_1 & \lambda_2 \\ &&& -\lambda_2 & a+\lambda_1 \end{bmatrix}$ and $B = \begin{bmatrix} u & p & q \\ & \mu_1 & \mu_2 && t \\ & -\mu_2 & \mu_1 & y & z \\ &&& u+\mu_1 & \mu_2 \\ &&& -\mu_2 & u+\mu_1 \end{bmatrix}$.
\end{enumerate}

The linearly nil-independent condition of $A$ and $B$ is as follows:
\begin{equation}\label{3.2}
	\begin{array}{|l |l |}
		\hline \text{Forms of $(A, B)$} & \text{Linearly nil-independent conditions} \\ \hline
		\hline \text{\eqref{form1} and \eqref{form3}} & \ra \begin{bmatrix} a & u \\ \lambda_1 & \mu_1 \\ \lambda_2&\mu_2 \end{bmatrix}=2 \\
		\hline \eqref{form2} & \ra \begin{bmatrix} a&u\\ \lambda&\mu\\ \end{bmatrix}=2 \\ \hline
	\end{array}
\end{equation}

To eliminate $\sigma_2, \sigma_4, \sigma_5$ we change $X' \co X + \alpha X_4 + \beta X_5$ 
and $Y' \co Y + \gamma X_4 + \delta X_5$. This transformation creates
\begin{equation}\label{3.3}
	\begin{array}{l l}
			[X',Y'] = & (\sigma_2 - t\beta + g\delta)X_2 +(\sigma_3 - y\alpha - z\beta + k\gamma + l\delta)X_3 \\
			& + \left[\sigma_4 - (u+v)\alpha - r\beta + (a+b)\gamma + f\delta\right]X_4 \\
			& + \left[\sigma_5 - x\alpha - (u+w)\beta + h\gamma + (a+c)\delta\right]X_5.
		\end{array}
\end{equation}
	
To destroy off-diagonal elements of $A$ and $B$, we use a basis changing which is an automorphism of $\g_{5,2}$ as follows
\begin{equation}\label{3.4}
	\begin{array}{l l l}
		\begin{bmatrix} X'_1 \\ X'_2 \\ X'_3 \\ X'_4 \\ X'_5 \end{bmatrix} = G_1\begin{bmatrix} X_1 \\ X_2 \\ X_3 \\ X_4 \\ X_5 \end{bmatrix},
		&& G_1 = \begin{bmatrix} 1& g_1 & g_2 \\ & 1 &&& g_3 \\ && 1 & g_4 & g_5 \\ &&& 1 \\ &&&& 1 \end{bmatrix} \in \Au(\g_{5,2}).
	\end{array}
\end{equation}
Transformation \eqref{3.4} will transform $A$ and $B$ into
\[
	\begin{array}{l l}
		G_1AG_1^{-1} = \begin{bmatrix} a & d' & e' & m^A & n^A \\ & b & f & s^A & g' \\ & h & c & k' & l' \\ &&& a+b & f \\ &&& h & a+c \end{bmatrix},
		& G_1BG_1^{-1} = \begin{bmatrix} u & p' & q' & m^B & n^B \\ & v & r & s^B & t' \\ & x & w & y' & z' \\ &&& u+v & r \\ &&& x & u+w \end{bmatrix},
	\end{array}		 
\]
in which
\begin{equation}\label{3.5}
	\begin{array}{l r r r r r r}
		d' = d & -(a-b)g_1 &  +hg_2 \\
		e' = e & +fg_1 & -(a-c)g_2 \\
		g' = g &&& +(a-b+c)g_3 && -fg_5 \\
		k' = k &&&& +(a+b-c)g_4 & +hg_5 \\
		l' = l &&& -hg_3 & +fg_4 & +ag_5 \\
		p' = p & -(u-v)g_1 & +xg_2 \\
		q' = q & +rg_1 & -(u-w)g_2 \\
		t' = t &&& +(u-v+w)g_3 && -rg_5 \\
		y' = y &&&& +(u+v-w)g_4 & +xg_5 \\
		z' = z &&& -xg_3 & +rg_4 & +ug_5
	\end{array}
\end{equation}
Afterwards, we destroy $m^A, n^A, s^A, m^B, n^B, s^B$ by changing $X' \co X - s^AX_1 + m^AX_2 + n^AX_3$ 
and $Y' \co Y -s^BX_1 + m^BX_2 + n^BX_3$. 
	
To normalize non-zero off-diagonal elements of $A$ and $B$, we also use an automorphism of $\g_{5,2}$ which is as follows
\begin{equation}\label{3.6}
	\begin{array}{l l l}
		\begin{bmatrix} X'_1 \\ X'_2 \\ X'_3 \\ X'_4 \\ X'_5 \end{bmatrix} = G_2\begin{bmatrix} X_1 \\ X_2 \\ X_3 \\ X_4 \\ X_5 \end{bmatrix},
		&& G_2 = \di (h_1, h_2, h_3, h_1h_2, h_1h_3) \in \Au(\g_{5,2}). \tag{$G_2$}
	\end{array}
\end{equation}
Transformation \eqref{3.6} will transform $A$ and $B$ into
\[
	\begin{array}{l l}
		G_2AG_2^{-1} = \begin{bmatrix} a & \frac{h_1d}{h_2} & \frac{h_1e}{h_3} \\ & b & \frac{h_2f}{h_3} && \frac{h_2g}{h_1h_3} \\ & \frac{h_3h}{h_2} & c & \frac{h_3k}{h_1h_2} & \frac{l}{h_1} \\ &&& a+b & \frac{h_2f}{h_3} \\ &&& \frac{h_3h}{h_2} & a+c \end{bmatrix}, &
		 G_2BG_2^{-1} = \begin{bmatrix} u & \frac{h_1p}{h_2} & \frac{h_1q}{h_3} \\ & v & \frac{h_2r}{h_3} && \frac{h_2t}{h_1h_3} \\ & \frac{h_3x}{h_2} & w & \frac{h_3y}{h_1h_2} & \frac{z}{h_1} \\ &&& u+v & \frac{h_2r}{h_3} \\ &&& \frac{h_3x}{h_2} & u+w \end{bmatrix}.
	\end{array}		 
\]

	\subsection{The structure matrices are of form \eqref{form1}}\label{ssec3.1}

	In this case, we have
	\[
		\begin{array}{l l}
			A = \begin{bmatrix} a & d & e \\ & \lambda_1 &&& g \\ && \lambda_2 & k & l \\ &&& a+\lambda_1 \\ &&&& a+\lambda_2 \end{bmatrix},
			& B = \begin{bmatrix}  u & p & q \\ & \mu_1 &&& t \\ && \mu_2 & y & z \\ &&& u+\mu_1 \\ &&&& u+\mu_2 \end{bmatrix}.
		\end{array}
	\]

	Transformation \eqref{3.4} transforms $A$ and $B$ by \eqref{3.5} which becomes
	\begin{equation}\label{3.5a}
		\begin{array}{l l r r r r r r}
			d' = d & -(a-\lambda_1)g_1 \\
			e' = e && -(a-\lambda_2)g_2 \\
			g' = g &&& +(a-\lambda_1+\lambda_2)g_3 \\
			k' = k &&&& +(a+\lambda_1-\lambda_2)g_4 \\
			l' = l &&&&& +ag_5 \\
			p' = p & -(u-\mu_1)g_1 \\
			q' = q && -(u-\mu_2)g_2 \\
			t' = t &&& +(u-\mu_1+\mu_2)g_3 \\
			y' = y &&&& +(u+\mu_1-\mu_2)g_4 \\
			z' = z &&&&& +ug_5
		\end{array}\tag{$3.5a$}
	\end{equation}
	and \eqref{3.1} becomes
	\begin{equation}\label{3.1a}
		\begin{cases}
			(a-\lambda_1)p=(u-\mu_1)d\\
			(a-\lambda_2)q=(u-\mu_2)e\\
			(a-\lambda_1+\lambda_2)t=(u-\mu_1+\mu_2)g\\
			(a+\lambda_1-\lambda_2)y=(u+\mu_1-\mu_2)k\\
			az=ul  
		\end{cases}\tag{$3.1a$}
	\end{equation}

	By \eqref{3.1a} and \eqref{3.5a} we can see that:
	\begin{equation}\label{*}
		\begin{array}{|c|c|c|}
			\hline \text{If} & \text{we choose} & \text{then} \\ \hline 
			\hline (a-\lambda_1)^2+(u-\mu_1)^2 \neq 0 & 
				g_1 \in \left\lbrace \frac{d}{a-\lambda_1}, \frac{p}{u-\mu_1} \right\rbrace & d' = p ' = 0 \\
			\hline (a-\lambda_2)^2+(u-\mu_2)^2 \neq 0 &
				g_2 \in \left\lbrace \frac{e}{a-\lambda_2}, \frac{q}{u-\mu_2} \right\rbrace & e' = q'=0 \\
			\hline (a-\lambda_1+\lambda_2)^2+(u-\mu_1+\mu_2)^2 \neq 0 & 
				g_3 \in \left\lbrace \frac{g}{\lambda_1+\lambda_2-a}, \frac{t}{\mu_1+\mu_2-u} \right\rbrace & g' = t' =0\\
			\hline (a+\lambda_1-\lambda_2)^2+(u+\mu_1-\mu_2)^2 \neq 0 &
					g_4 \in \left\lbrace \frac{k}{\lambda_2-\lambda_1-a}, \frac{y}{\mu_2-\mu_1-u} \right\rbrace & k' = y' =0 \\
			\hline a^2 + u^2 \neq 0 & g_5 \in \left\lbrace -\frac{l}{a}, -\frac{z}{u} \right\rbrace& l' = z' = 0 \\ \hline
		\end{array}\tag{$*$}
	\end{equation}

	According to \eqref{Mub-ineq2}, the center $\Z$ of $L$ satisfies $\dim \Z \leq 3$. However, Lie brackets of $\g_{5,2}$ show that there are only $X_4, X_5$ can belong to $\Z$. Moreover, if $X_4, X_5 \in \Z$ then $a=-\lambda_1=-\lambda_2$ and $u=-\mu_1=-\mu_2$ which conflict \eqref{3.2}. Therefore, we only have $\dim \Z = 1$ or $\dim \Z = 0$.
	
	\subsubsection{$\dim \Z = 1$}
	
	In this subcase, we have $\Z = \Sp \{X_4\}$ or $\Z = \Sp\{X_5\}$. However, if we interchange $X_2 \leftrightarrow X_3$ and $X_4 \leftrightarrow X_5$ then they will coincide. Therefore, without loss of generality, we can assume $\Z = \Sp \{X_4\}$, i.e.,
	\[
		\begin{array}{l l}
			a+\lambda_1 = 0=u+\mu_1, & (a+\lambda_2)^2+(u+\mu_2)^2 \neq 0.
		\end{array}
	\]
	
	Since \eqref{3.2} guarantees the left-hand side of \eqref{*}, its right-hand side is always valid. In other words, we can always transform the structure matrices into the following diagonal forms:
	\[
		\begin{array}{l l}
			A = \di (a, -a, \lambda_2, 0, a+\lambda_2), & B = \di(u, -u, \mu_2, 0, u+\mu_2).
		\end{array}
	\]
	
	First, we have $\sigma_2 = \sigma_3 = 0$ as $[A, B] = 0$. Next, we can choose appropriately $\beta, \delta$ in \eqref{3.3} to destroy $\sigma_5$, i.e., $\xy = \sigma_4X_4$. Since $a^2+u^2 \neq 0$, we can assume $a \neq 0$, otherwise, we interchange $X \leftrightarrow Y$. Thus, we normalize $a=1$ by scaling $X \to \frac{1}{a}X$ and then destroy $u$ by changing $Y' \co Y - uX$. Since $\mu_2 \neq 0$, we normalize $\mu_2 = 1$ by scaling $Y \to \frac{1}{\mu_2}Y$ and then destroy $\lambda_2$ by changing $X' \co X - \lambda_2Y$. It creates the following Lie algebras:
	\[
		\begin{array}{l l l l}
			L_1^\sigma \colon & A = \di (1, -1, 0, 0, 1), & B = \di (0, 0, 1, 0, 1), & \xy = \sigma X_4.
		\end{array}
	\]
	
	\begin{remark}
		Lie brackets of $L_1^\sigma$ can be easily read off due to their structure matrices. 
		Beyond the original ones of $\g_{5,2}$ and $\xy$, we have additionally
		\[
			\begin{array}{l l l l l}
				[X, X_1] = X_1, & [X, X_2] = -X_2, & [X, X_5] = X_5, & [Y, X_3] = X_3, & [Y, X_5] = X_5.
			\end{array}
		\]
		In our view, using structure matrices has an advantage that is a global view of the obtained Lie algebras' structures, such as decomposability or grouping Lie algebras for testing isomorphism (see Subsection \ref{ssec3.4} below), becomes more easier. Therefore, from now on, we use the structure matrices instead of Lie brackets.
	\end{remark}
	
	\subsubsection{$\dim \Z = 0$}
	
	In this subcase, $(a+\lambda_1)^2+(u+\mu_1)^2 \neq 0$ and $(a+\lambda_2)^2+ (u+\mu_2)^2 \neq 0$. Due to \eqref{*}, we can divide this subcase into two mutually-exclusive subcases as follows. Note that in two subcases below, we always have $[A, B] = 0$ which implies $\sigma_2 = \sigma_3 = 0$. Moreover, $\sigma_4, \sigma_5$ can always be eliminated by \eqref{3.3}, i.e., $\xy = 0$ in all two subcases.
	
	\paragraph{All of $d$, $e$, $g$, $k$, $l$, $p$, $q$, $t$, $y$, $z$ are zero}\label{221}
	
	It happens when all of $d$, $e$, $g$, $k$, $l$, $p$, $q$, $t$, $y$, $z$ are automatically zero or five inequalities on the left-hand side of \eqref{*} hold. This means that
	\[
		\begin{array}{l l}
			A = \di(a, \lambda_1, \lambda_2, a+\lambda_1, a+\lambda_2), &
			B = \di (u, \mu_1, \mu_2, u+\mu_1, u+\mu_2).
		\end{array}
	\]
	
	\begin{enumerate}[\bf A.]
		\item\label{A0} If $\lambda_1=\mu_1=0$ then the linearly nil-independent condition of $A$ and $B$ becomes $\ra \begin{bmatrix} a & u \\ \lambda_2 & \mu_2 \end{bmatrix} = 2$ which implies $a^2+u^2 \neq 0$ and $\lambda_2^2+\mu_2^2 \neq 0$. Without loss of generality, we can assume $a \neq 0$, otherwise, we interchange $X \leftrightarrow Y$. Thus, we normalize $a=1$ by scaling $X \to \frac{1}{a}X$ and then destroy $u$ by changing $Y' \co Y - uX$. Then, $\mu_2 \neq 0$, we normalize $\mu_2=1$ by scaling $Y \to \frac{1}{\mu_2}Y$ and then destroy $\lambda_2$ by changing $X' \co X - \lambda_2Y$.
 		
 		\item If $\lambda_1^2+\mu_1^2\neq0$ then we can assume $\lambda_2^2+\mu_2^2 \neq 0$ since on the contrary, we interchange $X_2 \leftrightarrow X_3$ and $X_4 \leftrightarrow X_5$ and return to \ref{A0}. If $a=u=0$ then we can assume $\lambda_1 \neq 0$. Thus, we normalize $\lambda_1=1$ by scaling $X \to \frac{1}{\lambda_1}X$ and then destroy $\mu_1$ by changing $Y' \co Y - \mu_1X$. Then, $\mu_2 \neq 0$, we normalize $\mu_2=1$ by scaling $Y \to \frac{1}{\mu_2}Y$ and destroy $\lambda_2$ by changing $X' \co X - \lambda_2Y$. If $a^2+u^2\neq0$ then we can assume $a \neq 0$, otherwise, we interchange $X \leftrightarrow Y$. Thus, we normalize $a=1$ by scaling $X \to \frac{1}{a}X$ and then destroy $u$ by changing $Y' \co Y- uX$. Then, $\mu_1^2+\mu_2^2 \neq 0$, we can assume $\mu_1 \neq 0$, otherwise, we interchange $X_2 \leftrightarrow X_3$ and $X_4 \leftrightarrow X_5$. Thus, we normalize $\mu_1=1$ by scaling $Y \to \frac{1}{\mu_1}Y$ and then destroy $\lambda_1$ by changing $X' \co X - \mu_1Y$.
	\end{enumerate}
	To sum up, we obtain the following Lie algebras:
	\[
		\begin{array}{l l l l}
			L_2 \colon & A = \di (1, 0, 0, 1, 1), & B = \di (0, 0, 1, 0, 1), \\ 
			L_3 \colon & A = \di (0, 1, 0, 1, 0), & B = \di (0, 0, 1, 0, 1), \\ 
			L^{ab}_4 \colon & A = \di (1, 0, a, 1, 1+a), & B = \di (0, 1, b, 1, b); & (a, b) \neq (0, -1).
		\end{array}
	\]
	
	\paragraph{There exists at least one of $d$, $e$, $g$, $k$, $l$, $p$, $q$, $t$, $y$, $z$ which is non-zero}
	
	It happens when $A$ or $B$ consists of non-zero off-diagonal elements $d$, $e$, $g$, $k$, $l$, $p$, $q$, $t$, $y$, $z$ and five inequalities on the left-hand side of \eqref{*} do not hold. However, \eqref{3.2} guarantees that there is at most one of them cannot hold. This means that $(A, B)$ can only contain at most one pair of non-zero off-diagonal elements which is $(d, p)$ or $(e, q)$ or $(g, t)$ or $(k, y)$ or $(l, z)$. Furthermore, if we interchange $X_2 \leftrightarrow X_3$ and $X_4 \leftrightarrow X_5$ then the pairs $(e, q)$ and $(k,y)$ will return to $(d, p)$ and $(g, t)$, respectively. Therefore, we have three situations as follows.
	\begin{enumerate}[\bf A.]
		\item \emph{$(A, B)$ contains the pair $(d, p)$}. We have $a = \lambda_1$ and $u=\mu_1$. First, we can assume $p \neq 0$, otherwise, we interchange $X \leftrightarrow Y$. Then we destroy $d$ by changing $X' \co X - \frac{d}{p}Y$. Note that \eqref{3.2} implies $a^2 + \lambda_2^2 \neq 0$. If $a = 0$, we normalize $\lambda_2 = 1$ by scaling $X \to \frac{1}{\lambda_2}X$ and then destroy $\mu_2$ by changing $Y' \co Y - \mu_2X$; afterwards, we normalize $u=p=1$ by $Y \to \frac{1}{u}Y$ and $G_2 = \di \left(\frac{u}{p}, 1, 1, \frac{u}{p}, \frac{u}{p}\right)$. If $a \neq 0$, we normalize $a = 1$ by scaling $X \to \frac{1}{a}X$ and then destroy $u$ by changing $Y' \co Y - uX$; afterwards, we normalize $\mu_2 = p = 1$ by scaling $Y \to \frac{1}{\mu_2}Y$ and $G_2 = \di \left(\frac{\lambda_2}{p}, 1, 1, \frac{\lambda_2}{p}, \frac{\lambda_2}{p}\right)$. To sum up, we obtain the following Lie algebras:
		\[
			\begin{array}{l l l}
				L_5 \colon & A = \di (0, 0, 1, 0, 1), & B = \begin{bmatrix} 1 & 1 \\ & 1 \\ && 0 \\ &&& 2 \\ &&&& 1 \end{bmatrix}, \\
				L^a_6 \colon & A = \di (1, 1, a, 2, 1+a), & B = \begin{bmatrix} 0 & 1 \\ & 0 \\ && 1 \\ &&& 0 \\ &&&& 1 \end{bmatrix}.
			\end{array}
		\]
		
		\item \emph{$(A, B)$ contains the pair $(g, t)$}. We have $\lambda_1=a+\lambda_2$ and $\mu_1=u+\mu_2$. First, we can assume $t \neq 0$ and then destroy $g$ by changing $X' \co X - \frac{g}{t}Y$. By a similar way as above, we obtain the following Lie algebras:
		\[
			\begin{array}{l l l}
				L_7 \colon & A = \di (0, 1, 1, 1, 1), & B = \begin{bmatrix} 1 \\ & 1 &&& 1 \\ && 0 \\ &&& 2 \\ &&&& 1 \end{bmatrix}, \\
				L^a_8 \colon & A = \di (1, 1+a, a, 2+a, 1+a), & B = \begin{bmatrix} 0 \\ & 1 &&& 1 \\ && 1 \\ &&& 1 \\ &&&& 1 \end{bmatrix}.
			\end{array}
		\]
	
		\item \emph{$(A, B)$ contains the pair $(l, z)$.} We have $a = u = 0$. First, we can assume $z \neq 0$ and then destroy $l$ by changing $X' \co X - \frac{l}{z}Y$. By a similar way as above, we obtain the following Lie algebras:
		\[
			\begin{array}{l l l}
				L_9 \colon & A = \di (0, 0, 1, 0, 1), & B = \begin{bmatrix} 0 \\ & 1 \\ && 0 & & 1 \\ &&& 1 \\ &&&& 0 \end{bmatrix}, \\
				L^a_{10} \colon & A = \di (0 , 1 , a, 1, a), & B = \begin{bmatrix} 0 \\ & 0 \\ && 1 && 1 \\ &&& 0 \\ &&&& 1 \end{bmatrix}.
			\end{array}
		\]
	\end{enumerate}
	
	\subsection{The structure matrices are of form \eqref{form2}}\label{ssec3.2}

	In this case, we first destroy $\mu_1$ by changing $Y' \co Y - \mu_1X$ to get
	\[
		\begin{array}{l l}
			A = \begin{bmatrix} a & d & e \\ & \lambda & 1 && g \\ && \lambda & k & l \\ &&& a+\lambda & 1 \\ &&&& a+\lambda \end{bmatrix}, &
			B = \begin{bmatrix} u & p & q \\ & \mu &&& t \\ && \mu & y & z \\ &&& u+\mu \\ &&&& u+\mu \end{bmatrix}.
		\end{array}
	\]
	
	Transformation \eqref{3.4} transforms $A$ and $B$ by \eqref{3.5} which becomes
	\begin{equation}\label{3.5b}
		\begin{array}{l l r r r r r r}
			d' & = & d & -(a-\lambda)g_1 \\
			e' & = & e & +g_1 & -(a-\lambda)g_2 \\
			g' & = & g &&& +ag_3 && -g_5 \\
			k' & = & k &&&& +ag_4 \\
			l' & = & l &&&& +g_4 & +ag_5 \\
			p' & = & p & -(u-\mu)g_1 \\
			q' & = & q && -(u-\mu)g_2 \\
			t' & = & t &&& +ug_3 \\
			y' & = & y &&&& +ug_4 \\
			z' & = & z &&&&& +ug_5
		\end{array}\tag{$3.5b$}
	\end{equation}
	and \eqref{3.1} becomes
	\begin{equation}\label{3.1b}
		\begin{cases}
			(a-\lambda)p=(u-\mu)d, \; (a-\lambda)q+d\mu_1=(u-\mu)e+p\\
			at+\mu_1l=ug+z, \; y=k\mu_1, \; ay=uk, \; az=ul
		\end{cases}\tag{$3.1b$}
	\end{equation}

	Now, \eqref{3.2} implies $(a-\lambda)^2+(u-\mu)^2 \neq 0$ and $a^2+u^2 \neq 0$. Taking account of \eqref{3.1b} and \eqref{3.5b}, we can choose $g_1 \in \left\lbrace \frac{d}{a-\lambda}, \frac{p}{u-\mu} \right\rbrace$, $g_2 \in \left\lbrace \frac{(a-\lambda)e+d}{(a-\lambda)^2}, \frac{q}{u-\mu} \right\rbrace$, $g_3 \in \left\lbrace \frac{k-al-a^2g}{a^3}, -\frac{t}{u} \right\rbrace$, $g_4 \in \left\lbrace -\frac{k}{a}, -\frac{y}{u} \right\rbrace$, $g_5 \in \left\lbrace \frac{k-al}{a^2}, -\frac{z}{u} \right\rbrace$ to destroy all $d$, $e$, $p$, $q$, $g$, $k$, $l$, $t$, $y$, $z$. Therefore, the structure matrices are transformed into
	\[
		\begin{array}{l l}
			A = \begin{bmatrix} a \\ & \lambda & 1 \\ && \lambda \\ &&& a+\lambda & 1 \\ &&&& a+\lambda \end{bmatrix}, &
			B = \di (u, \mu, \mu, u+\mu, u+\mu).
		\end{array}
	\]

	We have $\sigma_2 = \sigma_3 = 0$ as $[A, B] = 0$. Besides, we can choose $\alpha, \beta, \gamma, \delta$ in \eqref{3.3} to destroy $\sigma_4, \sigma_5$, i.e. $\xy = 0$. If $u=0$, we normalize $\mu=1$ by scaling $Y \to \frac{1}{\mu}Y$ and then destroy $\lambda$ by changing $X' \co X - \lambda Y$ and normalize $a=1$ by scaling $X \to \frac{1}{a}X$ and $G_2=(1,a,1,a,1)$; othewise, we normalize $u=1$ by scaling $Y \to \frac{1}{u}Y$ and then destroy $a$ by changing $X' \co X - aY$, and normalize $\lambda=1$ by scaling $X \to \frac{1}{\lambda}X$ and $G_2 = \di (1, \lambda, 1, \lambda, 1)$. We interchange $X \leftrightarrow Y$ to get a good look and obtain the following Lie algebras:
	\[
		\begin{array}{l l l}
			L_{11} \colon & A = \di (0, 1, 1, 1, 1), & B = \begin{bmatrix} 1 \\ & 0 & 1 \\ && 0 \\ 	&&& 1 & 1 \\ &&&& 1 \end{bmatrix}, \\
			L^a_{12} \colon & A = \di (1, a, a, 1+a, 1+a), & A = \begin{bmatrix} 0 \\ & 1 & 1 \\ && 1 \\ &&& 1 & 1 \\ &&&& 1 \end{bmatrix}.
		\end{array}
	\]
	
	\subsection{The structure matrices are of form \eqref{form3}}\label{ssec3.3}

	In this case, we first normalize $\lambda_2=1$ by scaling $X \to \frac{1}{\lambda_2}X$ and then destroy $\mu_2$ by changing $Y' \co Y -\mu_2X$ to get
	\[
		\begin{array}{l l}
			A = \begin{bmatrix} a & d & e \\ & \lambda_1 & 1 && g \\ & -1 & \lambda_1 & k & l \\ &&& a+\lambda_1 & 1 \\ &&& -1 & a+\lambda_1 \end{bmatrix},
			B = \begin{bmatrix} u & p & q \\ & \mu_1 &&& t \\ && \mu_1 & y & z \\ &&& u+\mu_1 \\ &&&& u+\mu_1 \end{bmatrix}.
		\end{array}
	\]
	
	Transformation \eqref{3.4} transforms $A$ and $B$ by \eqref{3.5} which becomes
	\begin{equation}\label{3.5c}
		\begin{array}{l l r r r r r r}
			d' & = & d & -(a-\lambda_1)g_1 & -g_2 \\
			e' & = & e & +g_1 & -(a-\lambda_1)g_2 \\
			g' & = & g &&& +ag_3 && -g_5 \\
			k' & = & k &&&& +ag_4 & -g_5 \\
			l' & = & l &&& +g_3 & +g_4 & +ag_5 \\
			p' & = & p & -(u-\mu_1)g_1 \\
			q' & = & q & & -(u-\mu_1)g_2 \\
			t' & = & t &&& +ug_3 \\
			y' & = & y &&&& +ug_4 \\
			z' & = & z &&&&& +ug_5
		\end{array}\tag{$3.5c$}
	\end{equation}
	and \eqref{3.1} becomes
	\begin{equation}\label{3.1c}
		\begin{cases}
			(a-\lambda_1)p = (u-\mu_1)d-q, \;
			(a-\lambda_1)q = (u-\mu_1)e+p \\
			y =-t, \;
			at = ug+2z, \;
			ay-z = uk, \; az=ul
		\end{cases}\tag{$3.1c$}
	\end{equation}
	
	By \eqref{3.1c} and \eqref{3.5c}, we take $g_1\in \left\lbrace \frac{d(a-\lambda_1)-e}{1+(a-\lambda_1)^2}, \frac{p}{u-\mu_1} \right\rbrace$ and $g_2 \in \left\lbrace \frac{d+e(a-\lambda_1)}{1+(a-\lambda_1)^2}, \frac{q}{u-\mu_1} \right\rbrace$ to destroy $d', e', p', q'$. Moreover, if $a^2 + u^2 \neq 0$, we take $g_3 \in \left\lbrace \frac{k-al-g(1+a^2)}{a(a^2+2)}, -\frac{t}{u} \right\rbrace$, $g_4 \in \left\lbrace \frac{g-al-k(1+a^2)}{a(a^2+2)}, -\frac{y}{u}  \right\rbrace$, $g_5 \in \left\lbrace \frac{g+k-al}{a^2+2}, -\frac{z}{u} \right\rbrace$ to further destroy $g', k', l', t', y', z'$. Therefore, we divide this case into two mutually-exclusive subcases as follows. Note that in two subcases below, we always have $[A, B] = 0$ and $\sigma_4, \sigma_5$ can always be eliminated by \eqref{3.3}, i.e., $\xy = 0$ in all two subcases.
	
	\subsubsection{All of $g$, $k$, $l$, $t$, $y$, $z$ are zero}\label{sssec3.3.1}
	
	It happens when $g$, $k$, $l$, $t$, $y$, $z$ are automatically zero or $a^2 + \mu^2 \neq 0$. This means that
	\[
		\begin{array}{l l}
			A = \begin{bmatrix} a \\ & \lambda_1 & 1 \\ & -1 & \lambda_1 \\ &&& a+\lambda_1 & 1 \\ &&& -1 & a+\lambda_1 \end{bmatrix}, &
			B = \di (u, \mu_1, \mu_1, u+\mu_1,  u+\mu_1).
		\end{array}
	\]
	
	Note that \eqref{3.2} implies $u^2 + \mu_1^2 \neq 0$. If $u = 0$, we normalize $\mu_1 =1$ by scaling $Y \to \frac{1}{\mu_1}Y$ and then destroy $\lambda_1$ by changing $X' \co X - \lambda_1Y$; otherwise, we normalize $u =1$ by scaling $Y \to \frac{1}{u}Y$ and then destroy $a$ by changing $X' \co X - aY$. We interchange $X \leftrightarrow Y$ to get a good look. It creates the following Lie algebras:
	\[
		\begin{array}{l l l}
			L_{13}^a \colon & A = \di (0, 1, 1, 1, 1), & B = \begin{bmatrix} a \\ & 0 &1 \\ & -1 & 0 \\ &&& a & 1 \\ &&& -1 & a \end{bmatrix}, \\ 
			L_{14}^{ab} \colon & A = \di (1, a, a, 1+a, 1+a), & B = \begin{bmatrix} 0 \\ & b & 1 \\ & -1 & b \\ &&& b & 1 \\ &&& -1 & b \end{bmatrix}.
		\end{array}
	\]
	
	\subsubsection{There exists at least one of $g$, $k$, $l$, $t$, $y$, $z$ which is non-zero}
	
	It happens when $A$ and $B$ consists of non-zero elements $g$, $k$, $l$, $t$, $y$, $z$ and $a = u =0$. In this subcases, \eqref{3.1c} gives $z=0$ and $y = -t$ and we have
	\[
		\begin{array}{l l}
			A = \begin{bmatrix} 0 \\ & \lambda_1 & 1 && g \\ & -1 & \lambda_1 & k & l \\ &&& \lambda_1 & 1 \\ &&& -1 & \lambda_1 \end{bmatrix}, &
			B = \begin{bmatrix} 0 \\ & \mu_1 &&& t \\ && \mu_1 & -t \\ &&& \mu_1 \\ &&&& \mu_1 \end{bmatrix}.
		\end{array}
	\]
	
	Now, we take $g_3 = g_4 = - \frac{l}{2}$ in \eqref{3.5c} to destroy $l$. Furthermore, we can destroy $g$ or $k$ by taking $g_5 = g$ or $g_5 = k$ in \eqref{3.5c}, respectively. However, if we change $X_2 \leftrightarrow X_3$ and $X_4 \leftrightarrow X_5$ then they will coincide. So, we take $g_5 = k$ to destroy $k$ and get
	\[
		\begin{array}{l l}
			A = \begin{bmatrix} 0 \\ & 0 & 1 && g \\ & -1 & 0 \\ &&& 0 & 1 \\ &&& -1 & 0 \end{bmatrix}, &
			B = \begin{bmatrix} 0 \\ & 1 &&& t \\ && 1 & -t \\ &&& 1 \\ &&&& 1 \end{bmatrix}.
		\end{array}
	\]
	
	Since $g^2+t^2 \neq 0$ to avoid subcase \ref{sssec3.3.1}, we normalize $t = 1$ by $G_2 = \di (t, 1, 1, t, t)$ if $g = 0$; otherwise, we normalize $g = 1$ by $G_2 = \di (g, 1, 1, g, g)$. It creates the following Lie algebras:
	\[
		\begin{array}{l l l}
			L_{15} \colon & A = \begin{bmatrix} 0 \\ & 0 & 1 \\ & -1 & 0 \\ &&& 0 & 1 \\ &&& -1 & 0 \end{bmatrix}, &
				B = \begin{bmatrix} 0 \\ & 1 &&& 1 \\ && 1 & -1 \\ &&& 1 \\ &&&& 1 \end{bmatrix}, \\
			L_{16}^a \colon & A = \begin{bmatrix} 0 \\ & 0 & 1 && 1 \\ & -1 & 0 \\ &&& 0 & 1 \\ &&& -1 & 0 \end{bmatrix}, &
				B = \begin{bmatrix} 0 \\ & 1 &&& a \\ && 1 & -a \\ &&& 1 \\ &&&& 1 \end{bmatrix}.
		\end{array}
	\]
		
	\subsection{Testing isomorphism}\label{ssec3.4}

	So far, we have done stage 1 in Subsection \ref{ssec2.1} for $\g_{5,2}$. 
	Subsections \ref{ssec3.1}, \ref{ssec3.2} and \ref{ssec3.3} show that we have constructed the list $\Ls_3$ 
	which consists of sixteen families of 7-dimensional indecomposable solvable $\F$-Lie algebras having nilradical $\g_{5,2}$. 
	However, $\Ls_3$ is not optimal since some families may be redundant. The goal of this subsection is to proceed stage 2, i.e., 
	to test isomorphism between the obtained Lie algebras by the technique pointed out in Subsection \ref{ssec2.2}. 
	There are two steps as follows.

	\begin{enumerate}
		\item The first step is to reduce $L_1^\sigma$, $L_4^{ab}$, $L_6^a$, $L_8^a$, $L_{10}^a$, $L_{12}^a$, $L_{13}^a$, $L_{14}^{ab}$ and $L_{16}^a$.
		\begin{enumerate}
			\item For $L_1^\sigma$, since $\di(1, \sigma, 1, \sigma, 1, 1, 1)$ is an isomorphism $L_1^{\sigma \neq 0} \cong L_1^1$, 
			we can reduce $\sigma$ to $\sigma \in \{0, 1\}$.
			\item For $L_4^{ab}$, the transformation
			\[
				\begin{bmatrix}
					-1 &&&&& -1 \\
					& 0 & -1 \\
					& -1 & 0 &&&& 1 \\
					&& 1 & 0 & 1 \\
					1 & 1& & 1 & 0 & 1 \\
					&&&&& 1 \\
					&&&&& a & b
				\end{bmatrix} \quad (b \neq 0)
			\]
			gives rise to an isomorphism $L_4^{ab} \cong L_4^{\left(-\frac{a}{b}\right)\frac{1}{b}}$. 
			This means that two pairs $(a, b)$ and $\left(-\frac{a}{b}, \frac{1}{b} \right)$ determine isomorphic algebras.
			\item $L_6^a$ and $L_8^a$ cannot be reduced since each parameter determines a unique Lie algebra, i.e., 
			$L_6^a \cong L_6^b$ and $L_8^a \cong L_8^b$ if and only if $a = b$.
			\item $L_{10}^a \cong L_{10}^\frac{1}{a}$ by
			\[
				\begin{bmatrix}
					\frac{1}{a} &&&&&& -\frac{1}{a} \\
					& 0 & a \\
					& a & 0 \\
					&&& 0 & 1 \\
					&&& 1 & 0 \\
					&&&&& a & 1 \\
					&&&&&& \frac{1}{a}
				\end{bmatrix} \quad (a \neq 0).
			\]
			\item $L_{12}^a$ cannot be reduced: $L_{12}^a \cong L_{12}^b$ if and only if $a = b$.
			\item We have $L_{13}^a \cong L_{13}^{-a}$ and $L_{14}^{ab} \cong L_{14}^{a(-b)}$ by $\di(1, -1, 1, -1, 1, 1, -1)$. 
			\item For $L_{16}^a$, we have $L_{16}^a \cong L_{16}^{-a}$ by $\di(1, -1, 1, -1, 1, -1, 1)$. 
		\end{enumerate}

		\item In the second step, we test isomorphism between Lie algebras in different families. To this end, we first group families into three groups with respect to their forms of structure matrices, i.e., three groups corresponding to Subsections \ref{ssec3.1}, \ref{ssec3.2} and \ref{ssec3.3}. Afterwards, these groups can be split into subgroups by using the dimensions of centers and further forms of structure matrices. By this way, we need to test isomorphism of families in the inside of eight groups as follows:
		\[
			\begin{array}{l l l}
				\A_{1.1} \co \{L_1\} & \A_{1.2} \co \left\lbrace L_2, L_3, L_4^{ab}\right\rbrace, & \A_{1.3} \co \{L_5, L_6^a\}, \\
				\A_{1.4} \co\{L_7, L_8^a\}, & \A_{1.5} \co \{L_9, L_{10}^a\}, & \A_2 \co \{L_{11}, L_{12}^a\}, \\
				\A_{3.1} \co \left\lbrace L_{13}^a, L_{14}^{ab}\right\rbrace, & \A_{3.2} \co \{L_{15}, L_{16}^a\}.
			\end{array}
		\]
		For $\A_{1.1}$, the test does not arise. To save more testing times, we further use the characteristic series of Lie algebras, i.e., the derived series and the lower central series. After checking these series, we do not need to test isomorphism in $\A_{1.2}$ because of different dimensions of ideals in the derived series. Consequently, we just need to test isomorphism in six remaining groups. The result is that we cannot reduce $\Ls_3$, i.e., $\Ls_3$ consists of 16 families of $\F$-Lie algebras.
	\end{enumerate}
	
	\begin{remark}
		All 16 families above are also valid over $\C$, except for $L_{13}^a$, $L_{14}^{ab}$, $L_{15}$ and $L_{16}^a$. They disappear since the Jordan block $\begin{bmatrix} a & b \\ -b & a \end{bmatrix}$ in the structure matrices does not exist over $\C$.
	\end{remark}

\section{Two main theorems}\label{sec4}

In this section, we formulate the two main results of this paper. The first one is four lists $\Ls_1$, $\Ls_2$, $\Ls_3$ and $\Ls_4$ which consist of 7-dimensional real and complex indecomposable solvable Lie algebras with nilradicals $(\g_1)^2 \oplus \g_3$, $\g_1 \oplus \g_4$, $\g_{5,2}$ and $\g_{5,4}$, respectively. Detailed computations of $\Ls_3$ are given in Section \ref{sec3}. For $\Ls_1$, $\Ls_2$ and $\Ls_4$, the computations are absolutely similar to that of Section \ref{sec3}, however, they are quite long. For convenience, we sum up these lists in Theorem \ref{thm4.1}.

\begin{theorem}\label{thm4.1}
	Four lists $\Ls_1$, $\Ls_2$, $\Ls_3$ and $\Ls_4$ are as follows.
	\begin{enumerate}
		\item $\Ls_1$ contains 99 (resp., 57) families of real (resp., complex) Lie algebras which are in Table \ref{tab1}.
		\item $\Ls_2$ contains 12 families of real and complex Lie algebras which are in Table \ref{tab2}.
		\item $\Ls_3$ contains 16 (resp., 12) families of real (resp., complex) Lie algebras which are in Table \ref{tab3}.
		\item $\Ls_4$ contains precisely one real and complex Lie algebra which is in Table \ref{tab4}.
	\end{enumerate}
	Tables \ref{tab1}, \ref{tab2}, \ref{tab3} and \ref{tab4} are given in Appendix.
\end{theorem}

Theorem \ref{thm4.1} presents a classification of 7-dimensional indecomposable solvable extensions of 
$(\g_1)^2 \oplus \g_3$, $\g_1 \oplus \g_4$, $\g_{5,2}$ and $\g_{5,4}$. As mentioned in Section \ref{sec1}, 
finite-dimensional indecomposable solvable extensions of $(\g_1)^5$, $\g_{5,1}$, $\g_{5,3}$, $\g_{5,5}$ 
and $\g_{5,6}$ were investigated. For the sake of completeness, we sum up these results here.
\begin{enumerate}
	\item {\bf Nilradical $(\g_1)^5$.} Ndogmo and Winternitz \cite{NW94} presented a procedure to classify 
	all finite-dimensional solvable Lie algebras with abelian nilradical. By this procedure, we obtain 31 (resp., 23) 
	families of 7-dimensional real (resp., complex) indecomposable Lie algebras with 5-dimensional abelian nilradicals.
	
	\item {\bf Nilradical $\g_{5,1}$.} The nilradical $\g_{5,1}$ is the 5-dimensional Heisenberg Lie algebra $\h_5$. 
	Rubin and Winternitz \cite[Table A2]{RW93} presented a table which consists of 27 (resp., 8) families of 7-dimensional 
	real (resp., complex) indecomposable solvable Lie algebras with nilradical $\h_5$.
	
	\item {\bf Nilradical $\g_{5,3}$.} \v Snobl and Kar\'asek \cite{SK10} classified solvable extension of nilradical 
	$\n_{n,3}$ in which $\g_{5,3} \cong \n_{5,3}$. Due to \cite[Theorem 2]{SK10}, there is precisely one 7-dimensional 
	real and complex indecomposable solvable Lie algebra with nilradical $\n_{5,3}$ as follows:
	$A = \di (1, 0, 1, 0, 1)$, $B = \di (2, 2, 1, 1, 0)$, $\xy = 0$.
	By the procedure in Section \ref{sec2}, we obtain the same algebra.
	
	\item {\bf Nilradical $\g_{5,5}$.} \v Snobl and Winternitz \cite{SW05} classified solvable extension of nilradical $\n_{n,1}$ in which  $\g_{5,5} = \n_{5,1}$. Due to \cite[Theorem 3]{SW05}, there is precisely one 7-dimensional real and complex indecomposable solvable Lie algebra with nilradical $\n_{5,1}$ as follows:
	$A = \di (3, 2, 1, 0, 1)$, $B = \di (1, 1, 1, 1, 0)$, $\xy = 0$.
	By the procedure in Section \ref{sec2}, we also obtain the same algebra.
	
	\item {\bf Nilradical $\g_{5,6}$.} \v Snobl and Winternitz \cite{SW09} classified solvable extension of nilradical $\n_{n,2}$ in which $\g_{5,6} \cong \n_{5,2}$. Due to \cite[Theorem 1]{SW09}, all real and complex solvable extensions of $\n_{n,2}$ must be $(n+1)$-dimensional. In other words, there is no 7-dimensional real and complex indecomposable solvable Lie algebra with nilradical $\n_{5,2}$. By the procedure in Section \ref{sec2}, we also obtain the same result.
\end{enumerate}

Combining all above results with Theorem \ref{thm4.1}, our second main result is:

\begin{theorem}\label{thm4.2}
	There are 188 (resp., 115) families of 7-dimensional real (resp., complex) indecomposable sovable Lie algebras with 5-dimensional nilradicals. These amounts are distributed as follows:
	\begin{center}
		\begin{tabular}{c c c c c c c c c c c c}
			\hline Nilradicals & $(\g_1)^5$ & $(\g_1)^2 \oplus \g_3$ & $\g_1 \oplus \g_4$ & $\g_{5,1}$ & $\g_{5,2}$ & $\g_{5,3}$ & $\g_{5,4}$ & $\g_{5,5}$ & $\g_{5,6}$ \\ \hline
			Over $\R$ & 31 & 99 & 12 & 27 & 16 & 1 & 1 & 1 & 0 \\ 
				Over $\C$ & 23 & 57 & 12 & 8 & 12 & 1 & 1 & 1 & 0 \\ \hline 
		\end{tabular} 
	\end{center}
\end{theorem}

\section{Complete classification of $7$-dimensional solvable Lie algebras}\label{sec5}

As mentioned in Section \ref{sec1}, for a $7$-dimensional Lie algebras $L$, its nilradical $N(L)$ satisfies $\dim N(L) \in \{4, 5, 6, 7\}$. Hindeleh and Thompson \cite{HT08}, Parry \cite{Par07} and Gong \cite{Gon98} classified $7$-dimensional Lie algebras $L$ such that $\dim N(L) \in \{4, 6, 7\}$.
Here, we summarize these results.

\begin{itemize}
	\item Gong \cite{Gon98} classified 7-dimensional indecomposable nilpotent Lie algebras over $\R$ and algebraically closed fields in which there are 149 (resp., 125) families of real (resp., complex) Lie algebras. Precisely isomorphic conditions for families containing parameters are also included.
	
	\item Parry \cite{Par07} classified 7-dimensional real indecomposable solvable Lie algebras with 1-codimensional nilradicals in which there are 594 families of real Lie algebras. By performing a procedure as in Remark \ref{rem2.1}, we obtained 525 families of complex Lie algebras. As mentioned in Example \ref{ex2.2}, Parry's list should be refined more to avoid redundancy. 
	
	\item Hindeleh and Thompson \cite{HT08} classified 7-dimensional real and complex indecomposable solvable Lie algebras with 4-dimensional nilradicals. There are 8 (resp., 2) families of real (resp., complex) Lie algebras. This results also should be refined more, in particular, parameters' conditions to avoid decomposability and redundancy. For example, the condition of parameters of algebra $7.2(ab)$ in \cite[Section 6]{HT08} should be $ab \neq 0$ instead of $a^2 + b^2 \neq 0$ since both $7.2(0b)$ and $7.2(a0)$ are decomposable. Furthermore, by using our testing isomorphism procedure in Subsection \ref{ssec2.2}, we can see that two pairs $(a, b)$ and $(b, a)$ determine isomorphic Lie algebras. Hence, we can reduce parameters to $a \geq b$ and $ab \neq 0$ (over $\R$) or $|a| \geq |b| > 0$ (over $\C$). The other families can also be refined appropriately.
\end{itemize}

To sum up, we have the following theorem:

\begin{theorem}\label{thm5.1}
The class of 7-dimensional solvable Lie algebras consists of 939 and 767 families of real and complex Lie algebras, respectively.
\end{theorem}



\section*{Funding}
\addcontentsline{toc}{section}{Funding}

The fourth author was funded by Vingroup Joint Stock Company and was supported by the Domestic Master/PhD Scholarship Programme of Vingroup Innovation Foundation (VINIF), Vingroup Big Data Institute (VINBIGDATA) under grant VINIF.2020.TS.46, and by the project SPD2019.01.37.


\appendix

\section*{Appendix}\label{sec6}
\addcontentsline{toc}{section}{Appendix}

All of the four following tables will consist of three columns as follows.
\begin{itemize}
	\item Algebras in the first column denoted by $L$ exist both over $\C$ and $\R$, 
	while those denoted by $R$ indicate that they only exist over $\R$.
	
	\item The second column contains triples $(A, B, \xy)$ in which the disappearance of $\xy$ means that $\xy = 0$.
	For convenience, we denote respectively by $(a_1, \ldots, a_5)$, $E_{ij}$ and $S_{ab}$ the diagonal matrix 
	$\di (a_1, \ldots, a_5)$, the 5-square matrix whose only non-zero entry is 1 	in row $i$ and column $j$, 
	and the block $\begin{bmatrix} a & b \\ -b & a \end{bmatrix}$.
	
	\item The final column contains additional conditions of families with parameters
	in which the notation ``$\equiv$'' means that 	these parameters yield isomorphic algebras.
	If there is no condition then parameters are arbitrary, 	and if ``$\equiv$'' disappears then 
	parameters are optimal in the sense that different parameters give rise to non-isomorphic algebras. 
	\item Conventions: $\sigma, \sigma' \in \{0, 1\}$, $\epsilon \in \{0, \pm 1\}$ and $\delta = \pm 1$.
\end{itemize}
For instance, $L_{1.1}^\sigma$ having structure $(1, -1, 0, 0, 0)$, $(0, 0, 0, 0, 1)$, $\sigma X_3 + X_4$ 
is a family of complex and real Lie algebras such that
\[
	\begin{array}{l l l l l}
		A = \di (1, -1, 0, 0, 0), & B = \di (0, 0, 0, 0, 1), & \xy = \sigma X_3 + X_4,
	\end{array}
\]
in which $\sigma \in \{0,1\}$. 
Similarly, $R_{3.1}^{a\sigma}$ having structure $(a, -a, 0, S_{01})$, $(0, 0, 0, 1, 1)$, $\sigma X_3$ with 
$a \geq 0$, $(a, \sigma) \neq (0, 0)$ indicates a family of real Lie algebras such that
\[
	\begin{array}{l l l l l}
		A = \begin{bmatrix} a \\ & -a \\ && 0 \\ &&& 0 & 1 \\ &&& -1 & 0 \end{bmatrix}, & B = \di (0, 0, 0, 1, 1), & \xy = \sigma X_3,
	\end{array}
\]
in which $a \geq 0$, $\sigma \in \{0,1\}$ and $(a, \sigma) \neq (0, 0)$.

\subsection*{Table \ref{tab1}}
\addcontentsline{toc}{subsection}{Table \ref{tab1}}

\begin{longtable}{p{.055\textwidth} p{.475\textwidth} p{.375\textwidth}}
	\caption*{Table 1: Solvable Lie algebras with nilradical $(\g_1)^2 \oplus \g_3$}\label{tab1} \\
	\hline No. & $(A, B, \xy)$ & Notes \endfirsthead
	\caption*{Table 1 (continued)}\\
	\hline No. & $(A, B, \xy)$ & Notes \\ \hline \endhead \hline
	$L_{1.1}^\sigma$ & $(1, -1, 0, 0, 0)$, $(0, 0, 0, 0, 1)$, $\sigma X_3 + X_4$ &\\ \hline
	$L_{1.2}$ & $(0, 0, 0, 1, 0)$, $(0, 0, 0, 0, 1)$, $X_3$ &\\ \hline
	$L_{1.3}^{ab\sigma}$ & $(1, -1, 0, 0, a)$, $(0, 0, 0, 1, b)$, $\sigma X_3$ 
		& $(a, b), (a, \sigma), (b, \sigma) \neq (0, 0)$; $(a, b, \sigma) \equiv \left(\pm \frac{a}{b}, \frac{1}{b}, \sigma\right)$ \\ \hline
	$L_{1.4}$ & $(0, 0, 0, 0, 1)$, $(1, -1, 0, 1, 0) + E_{14}$, $X_3$ \\ \hline
	$L_{1.5}^{a\sigma}$ & $(1, -1, 0, 1, a)$, $(0, 0, 0, 0, 1) + E_{14}$, $\sigma X_3$ \\ \hline
	$L_{1.6}$ & $(0, 0, 0, 0, 1)$, $(1, -1, 0, 0, 0) + E_{43}$, $X_4$ \\ \hline
	$L_{1.7}^{a\sigma}$ & $(1, -1, 0, 0, a)$, $(0, 0, 0, 0, 1) + E_{43}$, $\sigma X_4$ & $(a, \sigma) \equiv (-a, \sigma)$ \\ \hline
	$L_{1.8}^{ab}$ & $\bal (0, a, a, 0, 1), (1, b, 1+b, 0, 0), X_4 \ea$ 
		& $(a, b) \neq (0, -1)$; $(a, b) \equiv \left(-\frac{a}{b}, \frac{1}{b} \right)$ \\ \hline
	$L_{1.9}$ & $(0, 0, 0, 0, 1)$, $(0, 1, 1, 0, 0) + E_{14}$, $X_4$ \\ \hline
	$L_{1.10}^{a\sigma}$ & $(0, 1, 1, 0, a)$, $(0, 0, 0, 0, 1) + E_{14}$, $\sigma X_4$ \\ \hline
	$L_{1.11}$ & $(0, 1, 1, 0, 0)$, $(1, 0, 1, 0, 1) + E_{15}$, $X_4$ \\ \hline
	$L_{1.12}^a$ & $(1, a, 1+a, 0,1)$, $(0, 1, 1, 0, 0) + E_{15}$, $X_4$ \\ \hline
	$L_{1.13}^a$ & $(1, a, 1+a, 0, 1+a)$, $(0, 1, 1, 0, 1) + E_{53}$, $X_4$ & $a \neq -1$; $a \equiv \frac{1}{a}$ \\ \hline
	$L_{1.14}^{abc}$ & $(0, a, a, b, 1)$, $(1, c, 1+c, 0, 0)$  
		& $a, b \neq 0$; $(a, b, 0) \equiv \left(\frac{a}{b}, \frac{1}{b}, 0\right)$, 
			$(a, b, c) \equiv \left(-\frac{a}{c}, b, \frac{1}{c}\right) \equiv \left(-\frac{a}{bc}, \frac{1}{b}, \frac{1}{c}\right)$ \\ \hline
	$L_{1.15}^{abcd}$ & $(a, b, a+b, 0, 1)$, $(c, d, c+d, 1, 0)$ 
		& $(a, b), (c, d), (a+b, c+d) \neq (0, 0)$; $(a, b, c, d) \equiv (b, a, d, c) \equiv (c, d, a, b)$ \\ \hline
	$L_{1.16}^{ab}$ & $(0, 1, 1, 0, a)$, $(1, 0, 1, 1, b) + E_{14}$ & $(a, b) \neq (0, 0)$ \\ \hline
	$L_{1.17}^{ab}$ & $(1, a, 1+a, 1, b)$, $(0, 0, 0, 0, 1) + E_{14}$ & $a \neq -1$ \\ \hline
	$L_{1.18}^{abc}$ & $(1, a, 1+a, 1, b)$, $(0, 1, 1, 0, c) + E_{14}$ & $(b, c) \neq (0, 0)$ \\ \hline
	$L_{1.19}^{ab}$ & $(1, a, 1+a, 1+a, 0)$, $(0, b, b, b, 1) + E_{43}$  
		& $(a, b) \neq (-1, 0)$; $(a, b) \equiv \left(\frac{1}{a}, -\frac{b}{a}\right)$ \\ \hline
	$L_{1.20}^{abc}$ & $(a, b, a+b, a+b, 1)$, $(1, c, 1+c, 1+c, 0) + E_{43}$ 
		& $(a, b) \neq (0, 0)$, $(a+b, c) \neq (0, -1)$; $(a, b, c) \equiv \left(b, a, \frac{1}{c}\right)$ \\ \hline
	$L_{1.21}$ & $(1, 0, 1, 1, 1) + E_{14}$, $(0, 1, 1, 0, 0) + E_{15}$ \\ \hline
	$L_{1.22}$ & $(0, 0, 0, 0, 1) + E_{14}$, $(1, 1, 2, 1, 0) + E_{24}$ \\ \hline
	$L_{1.23}^{ab}$ & $(0, 1, 1, 0, 1) + E_{14} + E_{25}$, $(1, 0, 1, 1, 0) + aE_{14} + bE_{25}$ 
		& $(a, b) \equiv \left(\frac{1}{b}, \frac{1}{a}\right)$ \\ \hline
	$L_{1.24}$ & $(0, 1, 1, 0, 1) + E_{25}$, $(1, 0, 1, 1, 0) + E_{14}$ \\ \hline
	$L_{1.25}$ & $(0, 1, 1, 0, 1) + E_{14}$, $(1, 0, 1, 1, 0) + E_{25}$ \\ \hline
	$L_{1.26}^a$ & $(0, 0, 0, 0, 1) + E_{43}$, $(1, 0, 1, 1, 0) + E_{14} + aE_{43}$, $-X_2$ \\ \hline
	$L_{1.27}^a$ & $(0, 0, 0, 0, 1) + E_{14}$, $(1, 0, 1, 1, 0) + aE_{14} + E_{43}$, $X_2$ \\ \hline
	$L_{1.28}^{ab}$ & $(0, 0, 0, 0, 1) + E_{14} + E_{43}$, $(1, 0, 1, 1, 0) + aE_{14} + bE_{43}$, $(b-a)X_2$ \\ \hline
	$L_{1.29}$ & $(1, 0, 1, 1, 1)$, $(0, 1, 1, 0, 1) + E_{14} + E_{53}$ \\ \hline
	$L_{1.30}^a$ & $(1, 0, 1, 1, 1) + E_{53}$, $(0, 1, 1, 0, 1) + E_{14} + aE_{53}$ \\ \hline
	$L_{1.31}^a$ & $(1, 0, 1, 1, 1) + E_{14}$, $(0, 1, 1, 0, 1) + aE_{14} + E_{53}$ \\ \hline
	$L_{1.32}^{ab}$ & $(1, 0, 1, 1, 1) + E_{14} + E_{53}$, $(0, 1, 1, 0, 1) + aE_{14} + bE_{53}$ \\\hline
	$L_{1.33}$ & $(1, 0, 1, 1, 1) + E_{43}$, $(0, 1, 1, 1, 1) + E_{53}$ \\ \hline
	$L_{2.1}^\sigma$ & $(1, -1, 0, 0, 0) + E_{45}$, $(0, 0, 0, 1, 1)$, $\sigma X_3$ \\ \hline
	$L_{2.2}^{a\sigma}$ & $(0, 0, 0, 1, 1) + E_{45}$, $(1, -1, 0, a, a)$, $\sigma X_3$
		& $(a, \sigma) \neq (0, 0)$; $(a, \sigma) \equiv (-a, \sigma)$ \\ \hline
	$L_{2.3}^{a\sigma}$ & $(0, 1, 1, 0, 0) + E_{45}$, $(1, a, 1+a, 0, 0)$, $\sigma X_4$ 
		& $(a, \sigma) \equiv \left(\frac{1}{a}, \sigma\right)$ \\ \hline
	$L_{2.4}^{ab}$ & $(0, a, a, 1, 1) + E_{45}$, $(1, b, 1+b, 0, 0)$ 
		& $a \neq 0$; $(a, b) \equiv \left(-\frac{a}{b}, \frac{1}{b}\right)$ \\ \hline
	$L_{2.5}^{abc}$ & $(1, a, 1+a, 0, 0) + E_{45}$, $(b, c, b+c, 1,1)$ 
		& $(a, b+c) \neq (-1, 0)$; $(a, b, c) \equiv \left(\frac{1}{a}, c, b\right)$ \\ \hline
	$L_{2.6}^\sigma$ & $(1, 0, 1, 1, 1) + E_{45}$, $(0, 1, 1, 0, 0) + \sigma E_{15}$ \\ \hline
	$L_{2.7}^a$ & $(1, 0, 1, 1, 1) + E_{14} + E_{45}$, $(0, 1, 1, 0, 0) + aE_{15}$ \\ \hline 
	$L_{2.8}^{a\sigma}$ & $(0, 1, 1, 0, 0) + E_{45}$, $(1, a, 1+a, 1, 1) + \sigma E_{15}$ \\ \hline
	$L_{2.9}^{ab}$ & $(0, 1, 1, 0, 0) + E_{14} + E_{45}$, $(1, a, 1+a, 1, 1) + bE_{15}$ \\ \hline
	$L_{2.10}^{a\sigma}$ & $(0, 1, 1, 1, 1) + E_{45}$, $(1, a, 1+a, 1+a, 1+a) + \sigma E_{43}$ 
		& $(a, \sigma) \equiv \left(\frac{1}{a}, \sigma\right)$ \\ \hline
	$L_{2.11}^{ab}$ & $(0, 1, 1, 1, 1) + E_{45} + E_{53}$, $(1, a, 1+a, 1+a, 1+a) + bE_{43}$
		& $(a, b) \equiv \left(\frac{1}{a}, \frac{b}{a^3}\right)$ \\ \hline
	$R_{3.1}^{a\sigma}$ & $(a, -a, 0, S_{01})$, $(0, 0, 0, 1, 1)$, $\sigma X_3$ & $a \geq 0$, $(a, \sigma) \neq (0, 0)$ \\ \hline
	$R_{3.2}^{ab\sigma}$ & $(0, 0, 0, S_{a1})$, $(1, -1, 0, b, b)$, $\sigma X_3$
		& $a, b \geq 0$, $(b, \sigma) \neq (0, 0)$ \\ \hline 
	$R_{3.3}^{abc}$ & $(0, a, a, S_{b1})$, $(1, c, 1+c, 0, 0)$ & $a > 0$, $b \geq 0$ \\ \hline
	$R_{3.4}^{abcd}$ & $(a, b, a+b, S_{01})$, $(c, d, c+d, 1, 1)$ & $a, b \geq 0$, $(a+b, c+d) \neq (0, 0)$ \\ \hline
	$L_{4.1}^{a\epsilon}$ & $(0, 0, 0, 0, 1) + E_{12}$, $(0, 0, 0, 1, a)$, $\epsilon X_3$ 
		& $(a, \epsilon) \neq (0, 0)$; Over $\R$: $(a>0, \epsilon) \equiv \left(\frac{1}{a}, \epsilon\right)$, 
			$(a<0, \epsilon) \equiv \left(\frac{1}{a}, -\epsilon\right)$; Over $\C$: $\epsilon \in \{0, 1\}$, 
			$(a, \epsilon) \equiv \left(\frac{1}{a}, \epsilon\right)$ \\ \hline
	$L_{4.2}$ & $(1, 1, 2, 0, 0) + E_{12}$, $(0, 0, 0, 1, 0)$, $X_5$ \\ \hline
	$L_{4.3}^a$ & $(0, 0, 0, 1, 0) + E_{12}$, $(1, 1, 2, a, 0)$, $X_5$ & \\ \hline
	$L^{ab}_{4.4}$ & $(1, 1, 2, 0, a) + E_{12}$, $(0, 0, 0, 1, b)$ & $a, b \neq 0$; $(a, b) \equiv \left(-\frac{a}{b}, \frac{1}{b}\right)$ \\ \hline
	$L^{abc}_{4.5}$ & $(0, 0, 0, 1, a) + E_{12}$, $(1, 1, 2, b, c)$ & $(a, c) \neq (0, 0)$; $(a, b, c) \equiv \left(\frac{1}{a}, c, b\right)$ \\ \hline
	$L_{4.6}^a$ & $(1, 1, 2, 1, 0) + E_{12} + aE_{24}$, $(0, 0, 0, 0, 1) + E_{14}$ \\ \hline
	$L^a_{4.7}$ & $(0, 0, 0, 0, 1) + E_{12}$, $(1, 1, 2, 1, a) + E_{14}$ \\ \hline
	$L_{4.8}^{ab}$ & $(0, 0, 0, 0, 1) + E_{12} + E_{24}$, $(1, 1, 2, 1, a) + bE_{14}$ \\ \hline
	$L_{4.9}^a$ & $(1, 1, 2, 2, 0) + E_{12} + aE_{43}$, $(0, 0, 0, 0, 1) + E_{43}$ \\ \hline
	$L_{4.10}^a$ & $(0, 0, 0, 0, 1) + E_{12}$, $(1, 1, 2, 2, a) + E_{43}$ \\ \hline
	$L^{ab}_{4.11}$ & $(0, 0, 0, 0, 1) + E_{12} + E_{43}$, $(1, 1, 2, 2, a) + bE_{43}$ \\ \hline 
	$L_{5.1}^a$ & $(1, 1, 2, 0, 0) + E_{12} + E_{45}$, $(0, 0, 0, 1, 1) + aE_{45}$ \\ \hline
	$L_{5.2}^{ab}$ & $(0, 0, 0, 1, 1) + E_{12} + E_{45}$, $(1, 1, 2, a, a) + bE_{45}$ \\ \hline
	$R_{6.1}^{\sigma\sigma'}$ & $(0, 0, 0, S_{01})$, $(0, 0, 0, 1, 1) + \sigma E_{12}$, $\sigma'X_3$ & $(\sigma, \sigma') \neq (0, 0)$ \\ \hline
	$R_{6.2}^{a\epsilon}$ & $(0, 0, 0, S_{01}) + E_{12}$, $(0, 0, 0, 1, 1) + aE_{12}$, $\epsilon X_3$ & $a \geq 0$ \\ \hline
	$R_{6.3}^a$ & $(a, a, 2a, S_{01}) + E_{12}$, $(0, 0, 0, 1, 1)$ & $a > 0$ \\ \hline
	$R_{6.4}^{ab}$ & $(a, a, 2a, S_{b1})$, $(0, 0, 0, 1, 1) + E_{12}$ & $a > 0$, $b \geq 0$ \\ \hline
	$R_{6.5}^{ab}$ & $(0, 0, 0, S_{a1}) + E_{12}$, $(1, 1, 2, b, b)$ & $a \geq 0$ \\ \hline
	$R_{6.6}^{abc}$ & $(0, 0, 0, S_{b1}) + aE_{12}$, $(1, 1, 2, c, c) + E_{12}$ 
		& $a, b \geq 0$, $(a, c) \neq (0, 0)$ \\ \hline
	$R_{7.1}$ & $(S_{01}, 0, 0, 0)$, $(0, 0, 0, 0, 1)$, $X_4$ \\ \hline
	$R^\epsilon_{7.2}$ & $(S_{01}, 0, 0, 0)$, $(0, 0, 0, 0, 1) + E_{43}$, $\epsilon X_4$ \\ \hline
	$R^{a\sigma}_{7.3}$ & $(S_{01}, 0, 0, 0) + E_{43}$, $(0, 0, 0, 0, 1) + aE_{43}$, $\sigma X_4$
		& $a \geq 0$, $(a, \sigma) \neq (0, 0)$ \\ \hline 
	$R^{a\delta}_{7.4}$ & $(S_{01}, 0, a, 0)$, $(0, 0, 0, 0, 1)$, $\delta X_3$ & $a > 0$ \\ \hline
	$R^{ab\sigma}_{7.5}$ & $(S_{01}, 0, 0, a)$, $(0, 0, 0, 1, b)$, $\sigma X_3$
		& $(a, \sigma) \neq (0, 0)$, $b \neq 0$; $(a, b, \sigma) \equiv \left(\pm\frac{a}{b}, \frac{1}{b}, \sigma\right)$ \\ \hline
	$R^a_{7.6}$ & $(S_{a1}, 2a, 0, 0)$, $(0, 0, 0, 1, 0)$, $X_5$ & $a > 0$ \\ \hline
	$R^{ab}_{7.7}$ & $(S_{01}, 0, a, 0)$, $(1, 1, 2, b, 0)$, $X_5$ & $a \geq 0$, $(a, b) \neq (0, 0)$ \\ \hline
	$R_{7.8}$ & $(S_{01}, 0, 0, 0) + E_{43}$, $(1, 1, 2, 2, 0)$, $X_5$ \\ \hline
	$R^a_{7.9}$ & $(S_{01}, 0, 0, 0) + aE_{43}$, $(1, 1, 2, 2, 0) + E_{43}$, $X_5$ & $a \geq 0$ \\ \hline
	$R^{abc}_{7.10}$ & $(S_{a1}, 2a, 0, b)$, $(0, 0, 0, 1, c)$
		& $a, b, c \neq 0$; $(a, b, c) \equiv (-a, -b, c) \equiv \left(\pm a, \mp \frac{b}{c}, \frac{1}{c} \right)$ \\ \hline
	$R^{abcd}_{7.11}$ & $(S_{01}, 0, a, b)$, $(1, 1, 2, c, d)$
		& $(a, c), (b, d) \neq (0, 0)$; $(a, b, c, d) \equiv (\pm a, \pm b, c, d) \equiv (\pm b, \pm a, d, c)$ \\ \hline
	$R^{ab}_{7.12}$ & $(S_{a1}, 2a, 2a, b)$, $(0, 0, 0, 0, 1) + E_{43}$ & $a > 0$, $b \geq 0$ \\ \hline
	$R^{ab}_{7.13}$ & $(S_{01}, 0, 0, a)$, $(1, 1, 2, 2, b) + E_{43}$ & $a \geq 0$, $(a, b) \neq (0, 0)$ \\ \hline
	$R^{abc}_{7.14}$ & $(S_{01}, 0, 0, a) + E_{43}$, $(1, 1, 2, 2, b) + cE_{43}$ & $a, c \geq 0$, $(a, b) \neq (0, 0)$ \\ \hline
	$R_{7.15}$ & $(S_{01}, 0, 0, 0)$, $(1, 1, 2, 2, 2) + E_{43}$ \\ \hline
	$R_{7.16}$ & $(S_{01}, 0, 0, 0) + E_{43}$, $(1, 1, 2, 2, 2) + E_{53}$ \\ \hline
	$R^{\sigma\delta}_{8.1}$ & $(S_{01}, 0, 0, 0)$, $(0, 0, 0, 1, 1) + \sigma E_{45}$, $\delta X_3$ \\ \hline
	$R^{a\epsilon}_{8.2}$ & $(S_{01}, 0, 0, 0) + E_{45}$, $(0, 0, 0, 1, 1) + aE_{45}$, $\epsilon X_3$ & $a \geq 0$ \\ \hline
	$R^\sigma_{8.3}$ & $(S_{01}, 0, 0, 0)$, $(1, 1, 2, 0, 0) + E_{45}$, $\sigma X_4$ \\ \hline
	$R^{a\sigma}_{8.4}$ & $(S_{01}, 0, 0, 0) + E_{45}$, $(1, 1, 2, 0, 0) + aE_{45}$, $\sigma X_4$ & $a \geq 0$ \\ \hline
	$R^{ab}_{8.5}$ & $(S_{a1}, 2a, 0, 0) + E_{45}$, $(0, 0, 0, 1, 1) + bE_{45}$ & $a > 0$, $b \geq 0$ \\ \hline
	$R^{ab\sigma}_{8.6}$ & $(S_{01}, 0, a, a)$, $(1, 1, 2, b, b) + \sigma E_{45}$ & $a \geq 0$, $(a, b) \neq (0, 0)$ \\ \hline
	$R^{abc}_{8.7}$ & $(S_{01}, 0, a, a) + E_{45}$, $(1, 1, 2, b, b) + cE_{45}$ & $a, c \geq 0$, $(a, b) \neq (0, 0)$ \\ \hline
	$R^{ab}_{8.8}$ & $(S_{a1}, 2a, 2a, 2a) + E_{45}$, $(1, 1, 2, 2, 2) + E_{43} + bE_{45}$ 
		& 	$(a, 0) \equiv(c, 0)$; $(0, b) \equiv (0, -b)$; $	\left(a, \frac{1}{a}\right) \equiv \left(b, \frac{1}{b}\right)$;
			$(a, b) \equiv \left(a, \pm \frac{b}{1-ab}\right), \; ab \neq 1$;
			$(a, b) \equiv(c, d)$, $\frac{b}{d} = \pm \frac{1-ab}{1-cd}$, $acd \neq 0$, $ab, cd \neq 1$ \\ \hline
	$R^{abc}_{8.9}$ & $(S_{a1}, 2a, 2a, 2a) + E_{45} + E_{53}$, $(1, 1, 2, 2, 2) + bE_{43} + c(E_{45}+E_{53})$ & $a, c \geq 0$ \\ \hline
	$R^{a\sigma}_{9.1}$ & $(S_{01}, 0, a, a)$, $(0, 0, 0, S_{01})$, $\sigma X_3$ & $a \geq 0$, $(a, \sigma) \neq (0, 0)$ \\ \hline
	$R^{ab\epsilon}_{9.2}$ & $(S_{01}, 0, S_{0a})$, $(0, 0, 0, S_{1b})$, $\epsilon X_3$ 
		& $\bal a, b \geq 0, (a, \epsilon) \neq (0, 0) \ea$ \\ \hline
	$R^{ab}_{9.3}$ & $(S_{a1}, 2a, b, b)$, $(0, 0, 0, S_{01})$ & $a, b > 0$ \\ \hline
	$R^{abc}_{9.4}$ & $(S_{a1}, 2a, S_{0b})$, $(0, 0, 0, S_{1c})$ & $a, b > 0$, $c \geq 0$ \\ \hline
	$R^{abcd}_{9.5}$ & $(S_{01}, 0, S_{ab})$, $(1, 1, 2, S_{cd})$ & $a, d \geq 0$, $(a, b, c, d) \neq (0, 0, 0, 0)$ \\ \hline
	$R_{9.6}$ & $(S_{01}, 0, S_{01})$, $(1, 1, 2, 1, 1) + E_{14} + E_{25}$ \\ \hline
	$R^{abc}_{9.7}$ & $(S_{01}, 0, S_{01}) + E_{14} + aE_{15}$, $(1, 1, 2, 1, 1) + b(E_{14}+E_{25}) + c(E_{15} - e_{24})$
		& $(a, 0, 0) \equiv (d, 0, 0)$; $(a, b, c) \equiv (d, e, 0)$, $e \neq 0$, $bd - ae = \mp c$, $cd + e = \pm b$;
			$(a, b, c) \equiv (d, cd, c)$, $c \neq 0$, $(1-d^2)(b+ac) = 0$;
			$(a, b, c) \equiv (d, e, f)$, $f, df+e \neq 0$, $(e^2 + f^2)a + (f-de)b = \pm (df + e)c$, $(df + e)b \pm (f-de )c = \pm (e^2+f^2)$ \\ \hline
\end{longtable}

\subsection*{Table \ref{tab2}}
\addcontentsline{toc}{subsection}{Table \ref{tab2}}

\begin{longtable}{p{.055\textwidth} p{.475\textwidth} p{.375\textwidth}}
	\caption*{Table 2: Solvable Lie algebras with nilradical $\g_1 \oplus \g_4$}\label{tab2} \\
	\hline No. & $(A, B, \xy)$ & Notes \endfirsthead 
	\caption*{Table 2 (continued)} \\ 
	\hline No. & $(A, B, \xy)$ & Notes \\ \hline \endhead \hline
	\endlastfoot \hline
	$L_1$ & $(0, 0, 0, 0, 1)$, $(1, -2, -1, 0, 0)$, $X_4$ \\ \hline
	$L_2$ & $(0, 1, 1, 1, 0)$, $(1, 0, 1, 2, 0)$, $X_5$ \\ \hline
	$L_3^a$ & $(0, 1, 1, 1, 0)$, $(a, 0, a, 2a, 1)$ & $a \neq 0$ \\ \hline
	$L_4^{ab}$ & $(1, a, 1+a, 2+a, 0)$, $(0, b, b, b, 1)$ & $b \neq 0$; 
		$L_4^{ab \neq 0}$ should be added to Wang et al. \cite[Theorem 3]{WLD08} \\ \hline
	$L_5^a$ & $(1, 1, 2, 3, a)$, $(0, 0, 0, 0, 1) + E_{12}$ \\ \hline
	$L_6$ & $(0, 1, 1, 1, 0)$, $(1, 0, 1, 2, 1) + E_{15}$ \\  \hline
	$L_7^a$ & $(1, a, 1+a, 2+a, 1)$, $(0, 1, 1, 1, 0) + E_{15}$ \\ \hline
	$L_8^{a\delta}$ & $(0, 1, 1, 1, a)$, $(0, 0, 0, 0, 1) + \delta E_{24}$ & Over $\C$: $\delta = 1$ \\ \hline
	$L_9$ & $(0, 1, 1, 1, 1)$, $(1, 0, 1, 2, 0) + E_{25}$ & Should be added to Wang et al. \cite[Theorem 3]{WLD08} \\ \hline
	$L_{10}^a$ & $(1, a, 1+a, 2+a, a)$, $(0, 1, 1, 1, 1) + E_{25}$ & Should be added to Wang et al. \cite[Theorem 3]{WLD08} \\ \hline
	$L_{11}$ & $(0, 1, 1, 1, 1)$, $(1, 0, 1, 2, 2) + E_{54}$ \\ \hline
	$L_{12}^a$ & $(1, a, 1+a, 2+a, 2+a)$, $(0, 1, 1, 1, 1) + E_{54}$ \\ \hline
\end{longtable}

\subsection*{Table \ref{tab3}}
\addcontentsline{toc}{subsection}{Table \ref{tab3}}

\begin{longtable}{p{.055\textwidth} p{.475\textwidth} p{.375\textwidth}}
	\caption*{Table 3: Solvable Lie algebras with nilradical $\g_{5,2}$}\label{tab3} \\
	\hline No. & $(A, B, \xy)$ & Notes \endfirsthead
	\caption*{Table 3 (continued)} \\ 
	\hline No. & $(A, B, \xy)$ & Notes \\ \hline \endhead \hline
	$L_1^\sigma$ & $(1, -1, 0, 0, 1)$, $(0, 0, 1, 0, 1)$, $\sigma X_4$ & \\ \hline
	$L_2$ & $(1, 0, 0, 1, 1)$, $(0, 0, 1, 0, 1)$ \\ \hline
	$L_3$ & $(0, 1, 0, 1, 0)$, $(0, 0, 1, 0, 1)$ \\ \hline
	$L^{ab}_4$ & $(1, 0, a, 1, 1+a)$, $(0, 1, b, 1, b)$ & $(a, b) \neq (-1, 0)$; $(a, b) \equiv \left(-\frac{a}{b}, \frac{1}{b}\right)$ \\ \hline
	$L_5$ & $(0, 0, 1, 0, 1)$, $(1, 1, 0, 2, 1) + E_{12}$ \\ \hline
	$L^a_6$ & $(1, 1, a, 2, 1+a)$, $(0, 0, 1, 0, 1) + E_{12}$ \\ \hline
	$L_7$ & $(0, 1, 1, 1, 1)$, $(1, 1, 0, 2, 1) + E_{25}$ \\ \hline
	$L^a_8$ & $(1, 1+a, a, 2+a, 1+a)$, $(0, 1, 1, 1, 1) + E_{25}$ \\ \hline
	$L_9$ & $(0, 0, 1, 0, 1)$, $(0, 1, 0 , 1, 0) + E_{35}$ \\ \hline
	$L^a_{10}$ & $(0 , 1 , a, 1, a)$, $(0, 0, 1, 0, 1) + E_{35}$ & $a \equiv \frac{1}{a}$ \\ \hline
	$L_{11}$ & $(0, 1, 1, 1, 1)$, $(1, 0, 0, 1, 1) + E_{23} + E_{45}$ \\ \hline
	$L^a_{12}$ & $(1, a, a, 1+a, 1+a)$, $(0, 1, 1, 1, 1) + E_{23} + E_{45}$ \\ \hline
	$R_{13}^a$ & $(0, 1, 1, 1, 1)$, $(a, S_{01}, S_{a1})$ & $a \geq 0$ \\ \hline
	$R_{14}^{ab}$ & $(1, a, a, 1+a, 1+a)$, $(0, S_{b1}, S_{b1})$ & $b \geq 0$ \\ \hline
	$R_{15}$ & $(0, S_{01}, S_{01})$, $(0, 1, 1, 1, 1) +E_{25} - E_{34}$ \\ \hline
	$R_{16}^a$ & $(0, S_{01}, S_{01}) + E_{25}$, $(0, 1, 1, 1, 1) + a(E_{25} - E_{34})$ & $a \geq 0$ \\ \hline
\end{longtable}

\subsection*{Table \ref{tab4}}
\addcontentsline{toc}{subsection}{Table \ref{tab4}}

\begin{longtable}{p{.055\textwidth} p{.475\textwidth} p{.375\textwidth}}
	\caption*{Table 4: Solvable Lie algebras with nilradical $\g_{5,4}$}\label{tab4} \\
	\hline No. & $(A, B, \xy)$ & Notes \endfirsthead
	\hline No. & $(A, B, \xy)$ & Notes \endhead \hline
	$L_1$ & $(1, 0, 1, 2, 1)$, $(0, 1, 1, 1, 2)$ \\ \hline
\end{longtable}

\end{document}